\documentclass[12pt]{amsart} 
\usepackage{latexsym,amssymb, amsthm,amsfonts,amsmath,amssymb,amstext,amsthm,amssymb}
\usepackage[mathscr]{euscript}
\usepackage[abbrev]{amsrefs}

\usepackage[all]{xy}

\usepackage{hyperref}

 \hypersetup{colorlinks=true, urlcolor=green, citecolor=red, linkcolor=blue}

\newcommand{\coh}{{\rm Coh}} 
\newcommand{\csys}{{\rm CohSys}} 
\newcommand{\csysab}{{\rm CohSys}_\alpha^\beta} 
\newcommand{\cliff}{{\rm Cliff}} 
\newcommand{\eee}{\mathcal{E}}
\newcommand{\fab}{\mathscr{F}_\alpha^\beta}
\newcommand{\tab}{\mathscr{T}_\alpha^\beta}

\newcommand{\fff}{\mathcal{F}} 

\newcommand{\hh}{\mathcal{H}}
\newcommand{\hhm}{\mathcal{H}^{-1}}
\newcommand{\hhz}{\mathcal{H}^{0}}
\newcommand{\im}{{\rm im}}

\newcommand{\mua}{\mu_\alpha}

\newcommand{\mmp}{\mathfrak{m}_{P}} 
\newcommand{\nn}{\mathbb{N}}

\newcommand{\oo}{\mathcal{O}} 
\newcommand{\op}{\mathcal{O}_P} 

\newcommand{\ww}{\omega} 

\usepackage{amssymb}
\usepackage{pb-diagram}
\usepackage{enumitem} 
\usepackage[dvipsnames]{xcolor}
\usepackage{tikz-cd}
\usepackage{graphicx}
\usepackage{amsmath}
\numberwithin{equation}{section}
\usepackage[all]{xy}

\newtheorem{mthm}{Main Theorem}
\newtheorem{Theorem}{Theorem}[section]
\newtheorem{Proposition}[Theorem]{Proposition}
\newtheorem{Lemma}[Theorem]{Lemma}
\newtheorem{Corollary}[Theorem]{Corollary}
\newtheorem{Definition}[Theorem]{Definition}
\newtheorem{Remark}[Theorem]{Remark}

\numberwithin{equation}{section}
\newcommand{\rank}{\operatorname{rank}}

\newcommand{\R}{\mathbb{R}} \newcommand{\C}{\mathbb{C}}
\newcommand{\Z}{\mathbb{Z}} 
\newcommand{\Q}{\mathbb{Q}}
\newcommand{\oc}{\mathcal{O}_{C}}

\newcommand{\cata}{\mathscr{A}}

\newcommand{\cald}{\mathcal{D}}
\newcommand{\catd}{\mathscr{D}}
\newcommand{\cale}{\mathcal{E}}
\newcommand{\calf}{\mathcal{F}}
\newcommand{\calg}{\mathcal{G}}

\newcommand{\calp}{\mathcal{P}}

\newcommand{\lra}{\longrightarrow}
\newcommand{\stab}{\operatorname{Stab}}
\newcommand{\slice}{\operatorname{Slice}}
\newcommand{\prestab}{\operatorname{Stab}^{\rm pre}}
\newcommand{\stablf}{\operatorname{Stab}^{\rm lf}}
\newcommand{\Hom}{\operatorname{Hom}}

\DeclareMathOperator{\rk}{{rk}}

\DeclareMathOperator{\coker}{coker}
\DeclareMathOperator{\Min}{Min}

\begin{document}

\title[]{Stability conditions for coherent systems on Integral Curves}

\author{Marcos Jardim}
\address{Universidade Estadual de Campinas (UNICAMP) \\ Instituto de Matemática, Estatísitica e Computação Científica (IMECC) \\ Departamento de Matem\'atica \\
Rua S\'ergio Buarque de Holanda, 651\ \ 13083-970 Campinas-SP, Brazil.}
\email{jardim@unicamp.br}

\author{Leonardo Roa-Leguizam\'on}
\address{Universidade Estadual de Campinas (UNICAMP) \\ Instituto de Matemática, Estatísitica e Computação Científica (IMECC) \\ Departamento de Matem\'atica \\
Rua S\'ergio Buarque de Holanda, 651\ \ 13083-970 Campinas-SP, Brazil.} \email{leoroale@unicamp.br}

\author{Renato Vidal Martins}
\address{Departamento de Matemática, ICEx, UFMG Av. Ant\'onio Carlos 6627, 30123-970 Belo 
Horizonte MG, Brazil.}
\email{renato@mat.ufmg.br}
\thanks{}



\date{\today}

\baselineskip=16pt

\begin{abstract}
We present stability conditions for the category of coherent systems on an integral curve. We define a three-parameter family of pre-stability conditions in its derived category using tilting, and we then investigate when these conditions qualify as true stability conditions. Additionally, we examine the semistability of specific objects under these conditions, namely: torsion, free, and complete tilted systems, without relying on the support property.
\end{abstract}

\maketitle

\tableofcontents

\section{Introduction}

The study of vector bundles and their moduli space on algebraic curves is a central theme in algebraic geometry, with deep connections to mathematical physics, integrable systems, differential geometry, and number theory (see, for instance, \cites{BGMN2, Lu, Ty}). A crucial problem is to determine the geometry of the moduli space in terms of the existence and structure of its subvarieties.  One of the subvarieties that has been of great interest is Brill--Noether subvarieties.  A Brill--Noether locus is a subset of the moduli space whose points correspond to bundles having at least $k$ independent global sections. The primary goal of Brill-Noether theory is to study these subsets, specifically addressing questions related to non-emptiness, connectedness, irreducibility, dimension, singularities, and topological and geometric structures. For line bundles on curves (Classical Brill-Noether theory), many of these questions have been answered \cites{Ke,KL72,KL74,FL,GH}, but much less is known about higher rank vector bundles, or on singular curves  \cites{CR, BDG}.   A first approach to studying Brill--Noether subvarieties is by determining an upper bound for the dimension of the space of sections. Clifford's Theorem establishes that for any semistable vector bundle $E$ on a smooth curve of rank $n$ and degree $d$ with $0 \leq d/n \leq 2g-2$, $h^0(E)= {\rm {dim} }\, H^0(E)\leq d/2+n$ \cite{BGN}. 

A fruitful generalization of a vector bundle is the concept of a \textit{coherent system}. A coherent system can be thought of as a vector bundle together with a distinguished subspace of its global sections (also called Brill--Noether pair \cite{NMo}). More formally, a \textit{coherent system} of type $(n,d,k)$ on $C$ is a triple $\mathcal{E}=(E,V,\varphi)$ consisting of a coherent sheaf $E$ of rank $n$ and degree $d$, a $k$-dimensional vector space $V$, and a linear map $\varphi :V \to H^0(E)$. In addition, a morphism between coherent systems $(E,V,\varphi)$ and $(E',V',\varphi')$ is given by a linear map $f:V \to V'$ and a sheaf map $h:E \to E'$ such that $\varphi'\circ f=H^0(h)\circ \varphi$. In particular, if the vector bundle is stable, then the existence of such objects is equivalent to the non-emptiness of higher rank Brill--Noether subvarieties. This connection motivates the study of coherent systems (see, for instance, \cite{BGMN1}).

The category of coherent systems on $C$ denoted by $\csys(C)$ is also an example of a \emph{comma category}, see \cite[Example 3.2]{JNO} for details. In particular, $\csys(C)$ is abelian and noetherian, see for instance \cites{JNO, He, NMo}.

This simple definition significantly enriches the geometry. While the moduli space of stable vector bundles itself is a rich object, the theory of coherent systems allows us to stratify and study it by imposing conditions on the existence of families of sections. They provide a powerful framework for investigating the Brill--Noether theory of vector bundles. The study of linear systems, that is, rank 1 coherent systems for which the morphism $\varphi$ is injective, also has a vast literature, and it is closely related to the geometry of the curve itself; the best general reference here is \cite{ACGH}.  For a systematic treatment of coherent systems on smooth curves, we refer the reader to \cites{BGPMN, He, LP} 

 A fundamental step, as with vector bundles, is to define a notion of stability to construct well-behaved moduli spaces of coherent systems. The standard approach is to define the \textit{slope} of a coherent system. This notion of stability depends on a real parameter $\alpha$. Let $\alpha \in \mathbb{R}$. For a coherent system $\eee$ of type $(n, d, k)$, the \textit{$\alpha$-slope} is defined as:
\[ \mu_\alpha(E, V) = \frac{d}{n} + \alpha \frac{k}{n}. \]
A coherent system $\eee$ is called \textit{$\mu_\alpha$-(semi)stable} if for every proper coherent subsystem $\mathcal{F}$, the following inequality holds:
\[ \mu_\alpha(E', V') <(\le)~ \mu_\alpha(E, V). \]

The moduli space for coherent systems on curves was constructed by King and Newstead \cite{NMo}, and by Raghavendra and Vishwanath \cite{RaVi}. Moreover, the moduli space of coherent systems presents a rich display of topological and geometric phenomena in  Variational Geometric Invariant Theory (VGIT) (see \cites{BGPMN, Th} for more details). For a recent account of the theory and general open questions, we refer the reader to \cites{Ne, New}.

While much of the literature on coherent systems is devoted to smooth curves, we turn our attention to \textit{integral curves}, by which we mean a complete, integral, one-dimensional scheme over an algebraically closed field. One of the goals of this paper is to generalize some of the well-known results in the literature about coherent systems on smooth curves to coherent systems on integral curves.

A more recent trend in the theory of moduli spaces, initiated by T. Bridgeland almost 20 years ago, is the upgrading of stability conditions on abelian categories (like the category of coherent sheaves or the category of representations of a quiver) to triangulated categories, especially the derived categories of sheaves or representations. Bridgeland's notion of stability conditions on triangulated categories, introduced in \cite{Bridgeland} and \cite{Bridgeland1}, provides a new set of tools to study moduli spaces of sheaves on smooth projective varieties. Such tools have been successfully applied by many authors first to the study of sheaves on surfaces, for example, \cites{ABCH,BM14a,Bridgeland1,Macri-Schmidt}, and more recently on threefolds \cites{BMT,BMS,JLMM,Li,PT}. However, there has been little advance in the study of stability conditions on categories of \textit{decorated sheaves}; we are only aware of \cite{RHR} and the very recent preprint \cite{FN}, whose content partially overlaps with the present paper. The second goal of this paper is to provide a construction of stability conditions on the category of coherent systems and present certain classes of stable objects.

Let us present the organization of this paper and its main results. Section 2 is dedicated to revising the notation, basic definitions, and results regarding stability conditions and tilting theory. The following two sections focus on results related to curves, which will serve as the foundation for constructing stability conditions.

In Section \ref{sec:bounds}, we generalize Clifford’s Theorem for vector bundles to torsion-free sheaves that are not necessarily semistable on integral curves, as outlined in Theorem \ref{thmclf2}. The proof heavily relies on the rank 1 case, which was established in \cites{EKS,KM}. We then explore the conditions under which equality holds, referencing \cite{Ne}, and introduce Clifford indices while recalling some well-known inequalities for low slopes obtained in \cites{BGN,Mer}. 

Section \ref{sec:cs-integral} provides a brief summary of the main results concerning coherent systems on integral curves. We present a bound on the dimension of a coherent system in Theorem \ref{cota} and a generalization of Clifford's Theorem for coherent systems in Theorem \ref{thm:clifford-syst}. The former result will be crucial for defining a pre-stability condition in the derived category of coherent systems.

We finally arrive to the construction of stability conditions on $D^b(\csys(C))$ in Section \ref{sec:stab-cs}, where we introduce a two-parameter family of stability conditions $\sigma_\alpha^\beta=(\csys(C),Z_\alpha^\beta)$ whose heart is simply $\csys(C)$ regarded as the heart associated to the standard t-structure on $D^b(\csys(C))$; here $(\alpha,\beta)\in\R_{\ge0}\times\R$. For this reason, $\sigma_\alpha^\beta$ are called \textit{standard stability conditions}.

Sections \ref{sec:stab-tilt} and \ref{sec:bg-ineq} are dedicated to the construction of another family of (pre-)stability conditions using the technique called \textit{tilting on a torsion pair} that was successfully used in the context of stability conditions on the derived category of coherent sheaves on surfaces. It is worth mentioning here that Section \ref{sec:stab-tilt} has partial overlap with \cite[Section 3]{FN}. Given $\alpha\ge0$, and any $\beta\ge0$, we consider the torsion pair
$$ \tab := \{\mathcal{E}\in \csys(C)\, |\, \text{$\mu_\alpha(\mathcal{G}) > \beta$ whenever $\mathcal{E} \twoheadrightarrow \mathcal{G}$}\} $$
$$ \fab := \{\mathcal{E}\in \csys(C)\, |\, \text{ $\mu_\alpha(\mathcal{F}) \leq \beta$  for every $0\neq\mathcal{F} \subset\mathcal{E}$}\} $$
and the tilting category $\csysab(C):=\langle\fab[1],\tab\rangle$. When $\beta\ge0$, and $\gamma>1$, we show that the homomorphism $Z_\alpha^{\beta,\gamma}:K_0(\csysab(C))\to\C$ given by
$$ Z_\alpha^{\beta,\gamma}(\mathcal{E})= (d+\gamma n-k) + \sqrt{-1}\,(d+ \alpha \, k - \beta \, n) $$
is a stability function (i.e., the imaginary part is nonnegative, while the real part is negative when the imaginary one vanishes) on $\csysab(C)$. Therefore, we consider the pairs $\tau_\alpha^{\beta,\gamma}=(\csysab(C),Z_\alpha^{\beta,\gamma})$, the set
$\mathsf{PS}:=\R_{\ge0}\times\R_{\ge0}\times\R_{>1}$, and 
the subset $\mathsf{S}\subset\mathsf{PS}$ such that $\beta\ne1$ and either $\gamma\geq(\beta^2+2\alpha\beta-\alpha)/(\alpha(\beta-1))$, if $\beta >1$, or $\gamma > 1+(\beta(1-\alpha))/2\alpha$, if $\beta<1$.
The main results of Sections \ref{sec:stab-tilt} and \ref{sec:bg-ineq} can be summarized as follows; see Theorems \ref{StabCond} and \ref{supproperty}, Corollary \ref{cor:map2stab} and Lemmas \ref{action2} and \ref{action3} below for more details.

\begin{mthm} \label{mthm:stab}
If $(\alpha,\beta,\gamma)\in\mathsf{PS}\cap\Q$, then $\tau_\alpha^{\beta,\gamma}$ is a locally finite pre-stability condition. Moreover, if $(\alpha,\beta,\gamma)\in\mathsf{S}$, then $\tau_\alpha^{\beta,\gamma}$ is a stability condition, and the map \linebreak $\mathsf{S}\to\stab(D^b(\csys(C)))$ defined by $(\alpha,\beta,\gamma)\mapsto \tau_\alpha^{\beta,\gamma}$ is continuous. In addition,
\begin{itemize}
\item[(i)] For any $(\alpha,\beta,\gamma)\in\mathsf{PS}$, and $\alpha'\ge0$, $\tau_\alpha^{\beta,\gamma}\not\in\widetilde{GL}^+(2,\R)\cdot\sigma_{\alpha'}$;
\item[(ii)] For any $(\alpha,\beta,\gamma)$ and $(\alpha',\beta',\gamma')$ in $\mathsf{PS}$, then 
$\tau_\alpha^{\beta,\gamma}\not\in\widetilde{GL}^+(2,\R)\cdot\sigma_\alpha$;
\item[(iii)] if $(\alpha-\alpha')\gamma-(\alpha'+1)\beta \neq -\beta'(\alpha+1)$, then
$\tau_\alpha^{\beta,\gamma}\not\in\widetilde{GL}^+(2,\R)\cdot\tau_{\alpha'}^{\beta',\gamma'}$.
\end{itemize}
\end{mthm}

In particular, the tilted stability conditions $\tau_\alpha^{\beta,\gamma}$ are not equivalent, via the $\widetilde{GL}^+(2,\R)$ on $\stab(D^b(\csys(C)))$, to the standard ones $\sigma_\alpha^\beta$. 

Sections \ref{sec:sst-torsion} and \ref{sec:sst-complete} are dedicated to discovering examples of $\mu_\alpha^{\beta,\gamma}$-semistable objects and determining some walls. In particular, we study stability properties in the $\gamma$-direction and proved that objects $\mu_\alpha^{\beta,\gamma}$-semistable for $\gamma\gg 0$ are $\mu_\alpha$-semistable in $\csys(C)$. 

We summarize the results we got in the following theorem. To state it, we introduce some additional notation. For $\eee=(E,V,\varphi)\in\csys(C)$, we write $\eee=(E,V)$ if $\varphi$ is injective and we set $|\eee|:=(E,H^0(E))$. We say $\eee$ is \emph{complete} if $\eee=|\eee|$. Finally, for any $P\in C$, we denote by $\oo_P$ its structure, viewed as a closed subscheme of $C$.

\begin{mthm} 
\label{mthm:stable-objs}
The following holds:
\begin{itemize}
    \item[(i)] $(\oo_P,0)$, $(0,V)$, and $|\oo_P|$ are $\mu_\alpha^{\beta,\gamma}$-stable for every $(\alpha,\beta,\gamma)\in\mathsf{PS}$,  every point $P\in C$, and every $V$ $1$-dimensional $k$-vector space.
\item[(ii)] Let $\alpha,  \gamma \in \mathbb{Q}$ . Let $\eee \in \csys(C)$ be  $\mu_\alpha$-stable of type $(n,d,k)$, with $n>0$ and $-\gamma n<d$ if $k=0$. Assume  $\mu_\alpha(\mathcal{E})< \beta$. Assume also that $\beta < \delta_\eee$ if $n\geq 2$, where $\delta_\eee=min\{\mu_\alpha(\mathcal{Q}) \, | \, \mathcal{E} \twoheadrightarrow \mathcal{Q}\neq \eee\}$. Then; 
     \begin{itemize}
         \item [(a)] If $\eee=(E,V)$ is not complete, then $\mathcal{E}[1]$ is $\mu_\alpha^{\beta,\gamma}$-stable  in $\csys_\alpha^\beta(C)$ for every $ \gamma \gg 0$  and $\mu_\alpha^{\beta,\gamma}$-unstable for $ \gamma < \gamma_0:=(\beta n-d(\alpha+1))/\alpha n$. Also, if $|\eee| \in \fab$, then $\gamma_0$ is a pseudo-wall for $(n,d,k)$. Moreover, if $|\eee|$ is $\mu_\alpha$-stable and $\beta <\delta_{|\eee|}$, then  $\gamma_0$ is an actual-wall for $(n,d,k)$
         \item [(b)] If $\eee$ is  complete, then $\mathcal{E}[1]$ is $\mu_\alpha^{\beta,\gamma}$-stable for every $\gamma$ .
     \end{itemize}
    \end{itemize}
\end{mthm}

The significance of item (i) lies in its determination that our conditions are \emph{geometric} (see \cite[Definition 2.1]{BM11}). The proof of the entire result, which involves all items, relies on nearly all the findings established in this work. 

We informally refer to any object in the tilted category $\csysab(C)$ as a \emph{tilted system}. Accordingly, we call the objects in $\tab$ \emph{tilted torsion systems} and those in $\fab[1]$ \emph{tilted free systems}. This terminology explains the titles of Sections \ref{sec:stab-tilt}, \ref{sec:sst-torsion}, and \ref{sec:sst-complete}.

Section \ref{sec:sst-torsion} is dedicated to tilted torsion systems, with the goal of proving item (i). We begin our discussion in a broader categorical context, addressing the existence of minimal objects, elementary transformations, and their Harder-Narasimhan filtrations. When applied to our specific case, this topic naturally connects to \emph{base points} of coherent systems, as introduced in Definition \ref{defpb}. However, this notion, while useful for our purposes, does not fully capture the concept of a base point in a linear system when $C$ is singular. We provide a brief explanation of this right after the term is defined. Then we characterize semistability of tilted minimal objects in Lemma \ref{prpozz1}, and prove (i) in Theorem \ref{geometric} and Corollary \ref{opcomp}. The result relies on a description of coherent systems reaching maximal dimension, proved a few sections earlier in Proposition \ref{prop:extremal}.

We study tilted free systems in Section \ref{sec:sst-complete}, focusing on the $\mu_\alpha^{\beta,\gamma}$-semistability of these objects. Our goal now is to prove item (ii) mentioned above. We begin by examining the stability properties in the \(\gamma\)-direction. To this end, we establish the concepts of \emph{numerical}, \emph{pseudo}, and \emph{actual} walls in Definition \ref{critical}, drawing on insights from \cite[Section ~4]{JM}. We then identify necessary conditions for an object to be \(\mu_\alpha^{\beta,\gamma}\)-semistable when $\gamma$ is sufficiently large. Such an object should be a tilted $\mu_\alpha$-semistable torsion or free system, as stated in Lemma \ref{biggamma}. In Lemma \ref{subtor}, we describe the proper subobjects of a tilted system that satisfy the conditions of item (ii). Finally, in Theorem \ref{stabpos}, we address (a) and (b). As aforementioned, the proof relies on results from almost all sections of the present work.

\

\noindent
\textbf{Acknowledgements.} M.J. is partially supported by the CNPQ grant number 305601/2022-9 and the FAPESP-ANR Thematic Project 2021/04065-6. L.R-L. is partially supported by  FAPESP post-doctoral grant 2024/02475-0. R.V.M. is partially supported by CNPq grant number 308950/2023-2 and FAPESP grant number 2024/15918-8 


\section{Stability Conditions} \label{sec:stabconds}

This section mainly discusses stability conditions on triangulated categories. We derive several statements that will be referenced throughout this text. We begin by recalling important definitions and results for stability conditions on abelian categories.


\subsection{Stability conditions on abelian categories} \label{sec:stab-ab}

Let $\mathscr{A}$ be an abelian category with Grothendieck group $K_0(\mathscr{A})$, and $Z:K_0(\mathscr{A}) \to \mathbb{C}$ an additive homomorphism. Following \cite[Section 4]{Macri-Schmidt}, $Z$ is a called a \emph{stability function} if, for all nonzero $E \in \mathscr{A}$, we have
\[
\Im(Z(E)) \geq 0 \ \ \ \text{and} \ \ \ \Re(Z(E))<0\ \, \,  \text{if}\ \, \, \Im(Z(E))=0
\]
The \emph{slope} of $E\in \mathscr{A}$ (with respect to $Z$) is  
\[
M(E) := 
\begin{cases}
   \displaystyle{-\frac{\Re(Z(E))}{\Im(Z(E))}} \,\,\,\,\,\, \text{if $\Im(Z(E))\neq 0$} \\
    \infty \,\,\,\,\,\,\,\,\,\,\,\,\,\,\,\,\,\,\,\,\,\,\,\,\,\,\,\, \, \text{if $\Im(Z(E)) = 0$}
\end{cases}\]
The above definition suggests the following terms and notation: $D(E):=-\Re(Z(E))$ is the \emph{generalized degree} of $E$, while $R(E):=\Im(Z(E))$ is the \emph{generalized rank} of $E$.

A nonzero object $E\in \mathscr{A}$ is called \emph{stable} (respectively, \emph{semistable}) if, for all proper non trivial subobject $F \subset E$, we have $M(F)<M(E)$  (respectively, $\leq$). The \emph{phase} of $E$ is
\[\phi(E)= \frac{\rm arg(Z(E))}{\pi} \in (0,1].\]
It is easily seen that $E$ is stable (respectively, semistable) if, for all proper nontrivial subobjects $F\subset E$, we have $\phi(F)<\phi(E)$ (respectively, $\leq$).

\begin{Definition} \label{dfnstc}
A pair $(\mathscr{A},Z)$ is called a \emph{stability condition} if 
\begin{itemize}
\item[(i)] $Z$ is a stability function, and
\item[(ii)] every nonzero $E\in \mathscr{A}$ has a \emph{Harder--Narasimhan filtration}, that is, 
$$ 0=E_0 \subset E_1 \subset \ldots \subset E_n=E $$ 
with $E_i\in\mathscr{A}$, such that $E^i :=E_i/E_{i-1}$ is semistable for $1\leq i\leq n$, and
$$ M(E^1) > M(E^2) > \ldots > M(E^n) . $$
In addition. we define $M_{\rm max}(E):=M(E^1)$ and $M_{\rm min}(E):= M(E^n)$.
\end{itemize}
\end{Definition}

Finding sufficient conditions for the existence of Harder--Narasimhan filtrations in an abelian category is a crucial tool. The following criterion is stated in \cite[Prop.~4.9]{Macri-Schmidt}.

\begin{Proposition} \label{CHNF}
Let  $Z:K_0(\cata) \to \mathbb{C}$ be a stability function. Assume that
\begin{itemize}
\item [(i)] $\cata$ is Noetherian, and 
\item [(ii)] the image of $\Im(Z)$ is discrete in $\mathbb{R}$.
\end{itemize}
Then every nonzero $E\in\cata$ admits a Harder–-Narasimhan filtration with respect to $Z$.
\end{Proposition}

 The following useful result is proved in \cite[Lemma 4.5]{Macri-Schmidt}.
 
\begin{Proposition}\label{Homomorphism}
Let $Z:K_0(\cata) \to \mathbb{C}$ be a stability function. Let $A, B \in \cata$ be nonzero objects which are semistable with $M(A)>M(B)$. Then $\Hom(A,B)=0$.
\end{Proposition}

\subsection{Stability conditions on triangulated categories} \label{sec:stab-tri}

Bridgeland introduced stabi-\\lity conditions on triangulated categories in \cite[Definition 5.1]{Bridgeland} via \emph{slicings}. Later on, in \cite[Proposition 5.3]{Bridgeland}, he proves that this is equivalent to giving a bounded $t$-structure and a stability function on its heart satisfying the Harder--Narasimhan condition. To be more precise, we follow \cite[Section 5]{Macri-Schmidt}, introducing additional data, and applying it here to general triangulated categories.

Let $\mathscr{D}$ be a triangulated category, and fix a finite rank lattice $\Lambda$, a surjective group homomorphism $v:K_0(\mathscr{D})\to \Lambda$, and a norm $\left \| \cdot \right \|$ in $\Lambda$. Recall that if $\mathscr{A}$ is the heart of a bounded $t$-structure on $\mathscr{D}$, then $K_0(\mathscr{A})=K_0(\mathscr{D})$.  So $v$ is defined in $K_0(\mathscr{A})$ as well.

\begin{Definition}\label{BSC}
A \emph{pre-stability condition} on $\mathscr{D}$ is a pair $\sigma:=(\mathscr{A}, Z)$ where $\mathscr{A}$ is the heart of a bounded $t$-structure on $\mathscr{D}$, and $Z:K_0(\mathscr{A})\to\mathbb{C}$ is a stability function which factors through $v$, that is,
\begin{equation}\label{diadia}\begin{split}
\xymatrix{
K_0(\mathscr{A}) \ar[rr]^{Z} \ar[dr]^{v}& & \mathbb{C}\\
&\Lambda \ar[ur] & }   
\end{split} \end{equation}
such that $(\mathscr{A},Z)$ is a stability condition in the sense of Definition \ref{dfnstc}.

In addition, a \emph{stability condition} on $\mathscr{D}$ is a pre-stability condition satisfying the \emph{support property}:
\begin{equation} \label{equcsp}
S(\sigma):=\inf\left\{\frac{|Z(E)|}{\left \| v(E) \right \|}\ \bigg{|}\ 0 \neq E \in \mathscr{A} \,\,\text{semistable} \right\} > 0
\end{equation}
\end{Definition}

An alternative formulation of \eqref{equcsp} is given by the following useful result \cite[Exercise 5.9]{Macri-Schmidt}.

\begin{Proposition}
\label{prposp}
Let $\sigma=(\mathscr{A}, Z)$ be a pre-stability condition. Then $\sigma$ satisfies the support property \eqref{equcsp} if and only if there is a symmetric bilinear form $Q$ on $\Lambda_{\mathbb{R}}:=\Lambda \otimes \mathbb{R}$ such that: 
\begin{itemize}
    \item [(i)] $Q(\eee):=Q(v(\mathcal{E}), v(\mathcal{E}))\geq 0$ for all semistable $\mathcal{E} \in \mathscr{A}$;
    \item [(ii)] $Q(w):=Q(w, w)< 0$ for all non zero $w\in \Lambda_{\mathbb{R}}$ with $Z_\Lambda(w)=0$.
\end{itemize}
\end{Proposition}

Let $\prestab(\mathscr{D})$ and $\stab(\mathscr{D})$ denote the sets of pre-stability and stability conditions on $\mathscr{D}$, respectively,  with respect to $(\Lambda,v)$; clearly $\stab(\mathscr{D})\subset \prestab(\mathscr{D})$.  


Take $(\mathscr{A},Z)\in\prestab(\mathscr{D})$.  For each $\phi \in (0,1]$, define the additive subcategory of $\mathscr{D}$ 
\[\mathcal{P}(\phi):= \{E \in \mathscr{A} \, |\, \text{$E$ is semistable of phase $\phi$}\} \cup \{0\}\]
and extend this definition for all $\phi \in \mathbb{R}$ by setting
$\mathcal{P}(\phi + 1) := \mathcal{P}(\phi)[1]$. Given an interval $I\subset\mathbb{R}$, set, as in \cite[p.~327]{Bridgeland}, $\mathcal{P}(I)$ to be the extension-closed subcategory of $\mathscr{D}$ generated by the subcategories $\mathcal{P}(\phi)$ for $\phi\in I$.

A pre-stability condition $(\mathscr{A},Z)$ is said to be \textit{locally-finite} if for any $\phi\in\R$, there is $\epsilon>0$ such that $\calp((\phi-\epsilon,\phi+\epsilon))$ is of finite length, cf. \cite[Definition 5.7]{Bridgeland}; let $\stablf(\cald)$ denote the set of locally-finite pre-stability conditions. Moreover, when $(\mathscr{A}, Z)$ is a stability condition, it is not hard to check that $(\mathcal{P},Z_{\Lambda})$ is a stability condition via \textit{slicings} in the sense of \cite[Definition 5.1]{Bridgeland}.

This approach endows $\stab(\mathscr{D})$ with a metric. For $\sigma=(\mathcal{P}, Z)$ and $\tau=(\mathcal{Q}, W)$ in $\stab(\mathscr{D})$, one sets $d(\sigma,\tau):={\rm sup}\{d_S(\mathcal{P},\mathcal{Q}),\left \| Z-W \right \|\}$ where $d_S$ is defined in terms of Harder--Narasimhan filtrations for slicings as in \cite[Section 5.3, p.~23]{Macri-Schmidt}. Bridgeland proved in \cite[Theorem 1.2]{Bridgeland} the following landmark result.

\begin{Theorem} \label{dimStab}
The natural map $\stab(\mathscr{D}) \to \Hom(\Lambda, \mathbb{C})$ sending $(\mathscr{A},Z)$ to $Z$ is a local homeomorphism. In particular, $\stab(\mathscr{D})$ is a complex manifold of dimension ${\rm rank}(\Lambda)$.
\end{Theorem}

This presentation via slicings is also used in \cite[Lemma 5.5.4]{Bayer} to get the following result.
\begin{Lemma} 
\label{lemcsp}
The function $\stab(\mathscr{D}) \to \mathbb{R}_{>0}$ assigning $\sigma$ to $S(\sigma)$ as in \eqref{equcsp} is continuos.
\end{Lemma}

More generally, the set $\slice(\cald)$ of all slicings fo $\cald$ can also be equipped with a topology \cite[Section 6]{Bridgeland}; since
$$ \stablf(\cald) \subset \prestab(\cald) \subset \Hom(\Lambda,\C)\times\slice(\cald) $$
we also obtain topologies on $\stablf(\cald)$ and $\prestab(\cald)$ and continuous maps 
$$ \stablf(\cald)\to\Hom(\Lambda,\C) ~~ {\rm  and} ~~ \prestab(\cald)\to\Hom(\Lambda,\C) . $$
In addition, we have continuous inclusions $\stab(\cald)\subset\stablf(\cald)\subset\prestab(\cald)$.

Next, within the same framework, we will describe the action of a group on $\stab(\mathscr{D})$. Consider the universal cover $\widetilde{GL}^+(2,\mathbb{R})$ of $GL^+(2,\mathbb{R})$, whose elements are pairs
$(T, f)$, where $T \in GL^+(2,\mathbb{R})$, $f : \mathbb{R} \to \mathbb{R}$ is an increasing function with $f(x + 1) = f(x) + 1$ for all $x \in \mathbb{R}$, and such that the maps induced by $T$ and $f$ on $S^1=\mathbb{R}/2\mathbb{Z}=\mathbb{R}^2-\{(0,0)\}/\mathbb{R}_{>0}$ coincide. The group $\widetilde{GL}^+(2,\mathbb{R})$ acts on $\prestab_{\Lambda}(\mathscr{D})$ by 
$(T,f)\cdot (\mathcal{P},Z_\lambda) := (\mathcal{P}',Z')$, where $Z':=T^{-1}\circ Z_{\lambda}$ and $\mathcal{P}'(\phi)=\mathcal{P}(f(\phi))$.

The action of $\widetilde{GL}^+(2,\mathbb{R})$ on $\prestab(\mathscr{D})$ preserves the semistable objects, but relabels their phases. In addition, it also preserves the subsets $\stablf(\mathscr{D})$ and $\prestab(\mathscr{D})$.


\subsection{Torsion pairs and tilting} \label{sec:tilt}

A widely used technique to construct examples of stability conditions is called \textit{tilting}, which provides a new heart of a bounded t-structure out of a known one.
The procedure is described below for derived categories.

Let $\mathscr{A}$ be an abelian category. A \emph{torsion pair} on $\mathscr{A}$ is a pair of full additive subcategories $(\mathscr{T},\mathscr{F})$ of $\mathscr{A}$ such that:
\begin{itemize}
\item [(i)] For any $T \in \mathscr{T}$ and $F \in \mathscr{F}$, $Hom(T,F)=0$.
\item [(ii)] For any $E \in \mathscr{A}$, there is an exact sequence 
\[ 0 \to T \to E \to F \to 0 \]
in $\mathscr{A}$ with $T \in \mathscr{T}$ and $F \in \mathscr{F}$.
\end{itemize}

Set $\mathscr{D}:=D^b(\mathscr{A})$. The \emph{tilting} of $\mathscr{A}$ with respect to a torsion pair $(\mathscr{T},\mathscr{F})$ is the smallest extension closed full subcategory of $\mathscr{D}$ containing $\mathscr{T}$ and $\mathscr{F}[1]$. It is denoted by $\mathscr{A}^{\divideontimes}:=\langle\mathscr{F}[1], \mathscr{T}\rangle \subset \mathscr{D}$.

The next result can be derived by adjusting \cite[Lemma 6.3, Exercise 6.4]{Macri-Schmidt}.

\begin{Proposition}\label{newheart}
$\mathscr{A}^{\divideontimes}$ is the heart of a bounded $t$-structure on $\mathscr{D}$.
\end{Proposition}

Now let $(\mathscr{A},Z)$ be a stability condition and $M$ be the slope associated with $Z$. For each $\beta\in\mathbb{R}$, we define the following pair of subcategories on $\mathscr{A}$.
\begin{align}
\label{equtbt}
\mathscr{T}^\beta(\mathscr{A},Z):=\mathscr{T}^\beta  &:=\{E\in \mathscr{A}\, |\, \text{$M(G) > \beta$ whenever $E \twoheadrightarrow G$}\}  \\
\label{equfbt} 
\mathscr{F}^\beta(\mathscr{A},Z):=\mathscr{F}^\beta  & :=\{E\in \mathscr{A}\, |\, \text{ $M(F) \leq \beta$  for every $0\neq F \subset E$}\}
\end{align}

\begin{Remark}\label{filtrations}
It follows from \cite[Section 6.2 and Exercise 6.7]{Macri-Schmidt} that the categories $\mathscr{T}^\beta$ and $\mathscr{F}^\beta$ can be defined as
\begin{align}
\label{equtbt1}  \mathscr{T}^\beta  &:=\{E\in\cata\, |\, \text{any semistable factor $F$ of $E$ satisfies $M(F) > \beta$}\}  \\
\label{equfbt1}    \mathscr{F}^\beta  & :=\{E\in \cata\, |\, \text{ any semistable factor $F$ of $E$ satisfies $M(F) \leq \beta$}\}
\end{align}
\end{Remark} 

For future reference, we state the following result, whose proof is immediate.

\begin{Lemma} \label{prpelm}
Take $E\in\mathscr{A}$.
\begin{itemize}
\item[(i)] If $E$ is semistable and $M(E)\leq \beta$, then $E\in\mathscr{F}^\beta$;
\item[(ii)] if $E$ is semistable and $M(E) > \beta$, then $E\in\mathscr{T}^\beta$;
\item[(iii)] if $E\in\mathscr{F}^\beta$ has $M(E)=\beta$, then $E$ is semistable;
\item[(iv)] $\mathscr{F}^\beta$ is closed under subobjects;
\item[(iv)] $\mathscr{T}^\beta$ is closed under quotients.
\end{itemize}
\end{Lemma}


\section{Bounds on Global Sections of Sheaves on Integral Curves} \label{sec:bounds}

By a \emph{curve} we mean a complete, integral, one-dimensional scheme over an algebraically closed field. Throughout, $C$ stands for a given curve of arithmetic genus $g$, structure sheaf $\oc$, and dualizing sheaf $\ww$. Let $E$ be a coherent sheaf on $C$. Then $F:=E/{\rm Torsion}(E)$, as torsion-free, is locally free on the \emph{regular locus} $U$ of $C$ consisting of its non-singular points. Indeed, for every $P\in U$, the local ring $\mathcal{O}_P$ is a principal ideal domain, and apply \cite[Section 3.8]{J}. The \emph{rank} of $E$ is defined as the rank of $F|_U$. We set the \emph{degree} of $E$ as
\begin{equation} \label{equdgr}
\deg(E):=\chi(E)-\rk(E)\,\chi(\oc)
\end{equation}
Note that these definitions avoid the choice of a polarization on $C$; nonetheless, one can check that they agree with  \cite[Definitions 1.2.2, 1.2.11]{HL} or \cite[Definition 2.1]{KN}.

Let $\cata$ denote $\coh(C)$, the category of coherent sheaves on $C$. Consider the group homomorphism given by
$$ v: K_0(\cata) \to \Z^2 ~~,~~ [E] \mapsto (\rk(E),\deg(E)) ~~{\rm when}~~ E\in\cata, $$
which is easily seen to be surjective. We then define the additive function
$$ Z: K_0(C) \lra \C ~~,~~ Z([E]) = -\deg(E) + \sqrt{-1}\rk(E), $$
which factors through $\Z^2$; it is not hard to check that it is a stability function on $\cata$. The associated slope is just the usual Mumford slope
$$ \mu(E) = \dfrac{\deg(E)}{\rk(E)}. $$
Furthermore, as a consequence of Proposition \ref{CHNF} and Definition \ref{BSC}, the pair $(\cata,Z)$ defines a stability condition.

We are finally in a position to state the main result of this section.

\begin{Theorem} \label{thmclf1}
Let $E$ be a coherent sheaf on $C$ of rank $n$ and degree $d$. Then: 
\begin{itemize}
\item[(i)]  if $E$ is semistable and $d<0$, then $h^0(E)=0$;
\item[(ii)] if $\mu_{\rm min}(E)\ge0$, then
\begin{equation} \label{equhdn}
h^0(E)\leq d+n
\end{equation} 
\end{itemize}
\end{Theorem}

\begin{proof}
To prove (i), first assume $n=0$. Then \eqref{equdgr} yields $h^0(E)=d$ and we are done. So, for the remainder, consider $n>0$. Assume  $d<0$. If $h^0(E)>0$, then there is a morphism $\oc\to E$. But both $\oc$ and $E$ are semistable. Thus $\mu(\oc)\leq\mu(E)$ owing to \cite[Lemma 4.5]{Macri-Schmidt}. This yields $0\leq d/n<0$, a contradiction. So $h^0(E)=0$ when $d<0$.

For the proof of (ii), assume first that $E$ is semistable, so $\mu(E)=\mu_{\rm min}(E)\ge0$, implying that $d\ge0$. The exact sequence $0 \to E(-P) \to E \to \oo_P^{\oplus n} \to 0$
yields
\begin{equation} \label{equhzr}
h^0(E)\leq h^0(E(-P))+n
\end{equation} 
If $d<n$, then $\deg(E(-P))=d-n<0$. 
Since $E(-P)$ is semistable, the first item implies that $h^0(E(-P))=0$, hence $h^0(E)\leq n<d+n$. Thus (ii) holds for $d<n$. For $d\geq n$,
use induction on the degree: $\deg(E(-P))<d$ and $E(-P)$ is semistable. Then $h^0(E(-P))\leq(\deg(E)-n)+n$. So \eqref{equhzr} yields \eqref{equhdn}. 

When $E$ is not semistable, let $0\subset E_1\subset\cdots\subset E_l=E$ be its Harder--Narasimhan filtration, with factors $G_j:=E_j/E_{j-1}$. By assumption $\mu(G_j)\ge0$ so $d_j:=\deg(G_j)\ge0$; note that $\sum_j d_j=d$ and $\sum_j n_j=n$ where $n_j=\rk(G_j)$. We then have that
$$ h^0(E) \le \sum_j h^0(G_j) \le \sum_j (d_j+n_j) = d+n, $$
as desired.
\end{proof}

Our next result is a generalization of \cite[Theorem 2.1]{BGN}.

\begin{Theorem} \label{thmclf2}
If $E$ is coherent sheaf with $0\leq \mu_{\rm min}(E) \le \mu_{\rm max}(E) \leq 2g-2$, then
\begin{equation} \label{equcln} 
h^0(E)\leq \dfrac{\deg(E)}{2}+\rk(E)
\end{equation}
\end{Theorem}
\begin{proof}
We follow \cite[Theorem 2.1]{BGN}, making slight adjustments to the argument to suit our case, along with a few additional details. Note that the upper bound on $\mu_{\rm max}(E)$ implies that $E$ is torsion-free; otherwise, the maximal torsion subsheaf of $E$ is a semistable subsheaf with infinite slope.

Assume first that $E$ is semistable, so that $\mu(E)=\mu_{\rm min}(E)=\mu_{\rm max}(E)$.

Set $d=\deg(E)$ and $n=\rk(E)$. Note that if $h^1(E)=0$, then the definition of the degree \eqref{equdgr} yields $h^0(E)=d+n-ng = d/2+n-(ng-d/2)<d/2+n$. So \eqref{equcln} holds. So we may assume $h^0(E)\geq 1$ and $h^1(E)\geq 1$.

Now the proof goes by induction on the rank. If $n=1$, then apply \cite[Theorem A, p.~532]{EKS} or \cite[Lemma 3.1.(i)]{KM}. Now assume $n\geq 2$. Let $F\subset E$ be a torsion-free sheaf of maximal slope among all torsion-free subsheaves of $E$ of rank at most $n-1$. Then $F$ is semistable, as if $0\neq G\subset F$, then  $G$ is torsion-free, $G\subset E$, $\rank(G)\leq\rank(F)\leq n-1$, so $\mu(G)\leq\mu(F)$. Also, $E/F$ is semistable as well. Indeed, say $G/F\subset E/F$ where $F\subset G\subset E$ and consider the exact sequence $0\to F\to G\to G/F\to 0$. If $\rank(G)\leq n-1$, then $\mu(G)\leq\mu(F)$ as $G\subset E$ and $F$ is of maximal slope. Thus, the seesaw principle yields $\mu(G/F)\leq\mu(G)\leq\mu(F)$. But, similarly, as $E$ is semistable, from the exact sequence
\begin{equation}
\label{equfeq}
0\to F\to E \to E/F \to 0
\end{equation}
we get $\mu(F)\leq\mu(E)\leq\mu(E/F)$, so $\mu(G/F)\leq\mu(E/F)$. If $\rank(G)=n$, then $G/F\subset E/F$ are of the same rank, so  $\deg(G/F)\leq\deg(E/F)$, and hence $\mu(G/F)\leq\mu(E/F)$. Thus $E/F$ is semistable.

Now, by construction, $\mu(F)\geq 0$. Also, $\mu(F)\leq\mu(E)\leq 2g-2$. On the other hand, as $h^1(E)\geq 1$, there is a nonzero morphism $f\in{\rm Hom}(E,\ww)$. Consider the exact sequence
\begin{equation}
\label{equkei}
0\to \ker(f)\to E \to \im(f)\to 0.
\end{equation}
Note that $\im(f)$, as a nonzero coherent subsheaf of $\ww$, is torsion-free of rank 1. Thus $\rank(\ker(f))=n-1$. Also, $\deg(\im(f))\leq\deg(\ww)=2g-2$. Set $n':=\rank(E/F)$. So comparing slopes in \eqref{equfeq} and \eqref{equkei}, and using the fact that degrees are additive in exact sequences, yields
\begin{align*}
n'\mu(E/F) & =(n-1)\mu(\ker(f))+\deg(\im(f))-(n-n')\mu(F)\\
 & \leq (n-1)\mu(F)+2g-2-(n-n')\mu(F)=(n'-1)\mu(F)+2g-2 \\
 &\leq (n'-1)(2g-2)+2g-2=n'(2g-2)
\end{align*}
where the first inequality holds since $\mu(\ker(f))\leq\mu(F)$ because $\rank(\ker(f))=n-1$ and $F$ is of maximal slope. Thus $\mu(E/F)\leq2g-2$. Therefore, both $F$ and $E/F$ satisfy the same hypothesis of $E$, and are of smaller rank. Set $d':=\deg(E/F)$. Thus \eqref{equfeq} and induction yields $h^0(E)\leq h^0(F)+h^0(E/F)\leq (d-d')/2+(n-n')+d'/2+n'=d/2+n$.

When $E$ is not semistable, let $0\subset E_1\subset\cdots\subset E_l=E$ be its Harder--Narasimhan filtration, with torsion-free factors $E^j:=E_j/E_{j-1}$ of rank $n^j$ and degree $d^j$. By assumption, $0\le\mu(E^j)\le2g-2$ for each $j$, so $h^0(E^j)\le d_j/2+n_j$. We then have that
$$ h^0(E) \le \sum_j h^0(E^j) \le \sum_j \left(\frac{d^j}{2}+n^j\right) = \frac{d}{2}+n, $$
thus concluding the proof.
\end{proof}

We can now return to the characterization of the extremal cases in Theorem \ref{thmclf1}, which will be useful in Section \ref{sec:cs-integral} below.

\begin{Proposition}\label{charsheaf}
A semistable sheaf $E$ satisfies $h^0(E)=d+n$ if and only if either $n=0$, or $g=0$, or $E=\oc^{\oplus n}$.
\end{Proposition}

\begin{proof}
As previously noticed, equality in \eqref{equhdn} trivially holds if $n=0$. For $n\geq 1$, rewrite \eqref{equdgr} as $h^0(E)=d+(1-g)n+h^1(E)$. So equality holds in \eqref{equhdn} if and only if $h^1(E)=ng$.

If $g=0$, then $C=\mathbb{P}^1$ and any semistable vector bundle $E$ of rank $n$ and degree $d\geq 0$ on $C$ is of the form $E\cong \mathcal{O}_{\mathbb{P}^1}(a)^{\oplus n}$, with $na=d$. Hence $h^0(E)=\sum_{i=1}^n (a+1)=d+n$. So equality holds when $g=0$. 

Now assume $g>0$. Suppose $E$ is semistable, $h^1(E)=ng$, and $d\geq 0$. In particular, $h^1(E)>0$, so there is a nonzero morphism $E\to\ww_C$, and, as $E$ and $\ww_C$ are semistable, then $d/n\leq 2g-2$. Hence \eqref{equcln} holds for $E$. Therefore $h^0(E)\leq d/2+n\leq d+n=h^0(E)$, so $d=0$ and $h^0(E)=n$. As $E$ is semistable and $n>1$, then $E$ is torsion-free. As $h^0(E)=n$, by \cite[Proposition 3.1]{BGN}, there is an injection $\iota:\oc^{\oplus n}\hookrightarrow E$. But $E$ and $\oc^{\oplus n}$ are of the same degree zero and of the same rank, thus $E\simeq\oc^{\oplus n}$ since the $\coker(\iota)$ has zero rank and degree. Now note that $\oc^{\oplus n}$ is semistable whatever $C$ is, and this completes the proof.
\end{proof}

The characterization of equality in \eqref{equcln} requires more careful consideration. For example, the case when $n=1$ was addressed in \cite[Theorem A.(2), p.~532]{EKS}, but it involved a lengthy proof. This work introduced a new class of curves that satisfy Clifford's equality, distinct from hyperelliptic curves. Specifically, these are rational curves that possess a unique singularity whose conductor agrees with the maximal ideal. These curves are non-Gorenstein and were characterized in \cite[Theorem 5.10]{KM}. They are referred to as \emph{nearly normal} as they are the only curves whose canonical model is arithmetically normal.
 
For the case $n>1$, we may reframe the problem in terms of Clifford indices. Namely, Lange and Newstead introduced the \emph{Clifford index} of higher rank vector bundles on a smooth curve. As our $C$ is integral, we may replace locally free by torsion-free sheaves, as did, for instance, in \cite[Definition ~3.1]{LGMS} for $n=1$. So let $E$ be a torsion-free sheaf of rank $n$ and degree $d$ on $C$. In \cite[p.~169]{LN1}, the \emph{Clifford index} of $E$ is defined as follows
\begin{equation}
\label{equcle}
     \cliff(E)  := \frac{1}{n}(d - 2(h^{0}(E)-n)).
\end{equation}
So $\cliff(E)$ measures how far $E$ is from satisfying equality in \eqref{equcln}, as this holds if and only if $\cliff(E)=0$. If we want the Clifford index of a curve to measure how far it is from having such an $E$, then, accordingly, the \emph{$n$-Clifford index} of $C$ could be set as
$$
   \cliff_n(C)    :=  \min \left\{ \cliff\,(E) \,|\, E\ \text{is semistable of rank $n$ and $0\leq \mu(E)\leq 2g-2$} \right\}
$$
However, \cite[p.~169]{LN1}, where the invariant was coined, brings a more refined definition:
\begin{equation}
\label{equdcl}
\cliff_n(C)  :=  \min \left\{ \cliff\,(E) \,|\, E\ \text{is semistable of rank $n$},\ 0\leq \mu(E)\leq g-1, \, h^{0}(E) \geq 2n  \right\}.
\end{equation}
Indeed, rewrite \eqref{equcle} as
\begin{equation}
\label{equclg}
\cliff(E)=(g-1)-\frac{h^0(E)+h^1(E)}{n}
\end{equation}
and set $E^{\vee}_{\ww}:=\mathcal{H}om(E, \omega)$. So Serre's duality \cite[Theorem 7.6]{H}, along with \eqref{equclg}, easily yields $\cliff(E)=\cliff(E^{\vee}_\ww)$. It also yields $\mu(E)+\mu(E^{\vee}_\ww)=2g-2$. From this equality, we conclude that $E$ is semistable if and only if $E^{\vee}_\ww$ is, and, clearly, that $\mu(E)>g-1$ if and only if $\mu(E^{\vee}_\ww)\leq g-1$. So this explains the range $0\leq \mu(E)\leq g-1$.

Now rewrite \eqref{equcle} as
$$
\cliff(E)=\mu(E)-\frac{2h^0(E)}{n}+2.
$$
As $\cliff(E)\geq 0$, it follows that if $\mu(E)<2$, then $h^0(E) < 2n$. But if $\mu(E)< 2$, and $C$ is smooth, then $h^0(E)\leq (d-n)/g+n$, which is way stronger than \eqref{equcln}. Brambila--Grzegorczyk--Newstead proved this inequality in \cite[Theorem $\widetilde{\text{B}}$]{BGN} in the range $0\leq\mu(E)\leq 1$ and extended by Mercat in \cite[Théorème A-1, p.~11]{Mer} to the range $1<\mu(E)<2$. So this explains the range $h^0(E)\geq 2n$ in \eqref{equdcl}.  

In other words, \eqref{equdcl} may work for integral curves as well. Now, Mercat conjectured in \cite{Mer2} that $\cliff_n(C)=\cliff_1(C)$. In particular, when the statement holds, the characterization of curves with a torsion-free sheaf reaching the maximal number of independent global sections would reduce to the case $n=1$. For a recent account of this conjecture, we refer the reader to \cite{New}.


\section{Coherent Systems on Integral curves} \label{sec:cs-integral}

After discussing relevant definitions and results concerning coherent sheaves on an integral curve, we can now introduce the main subject of this work.

\begin{Definition}\label{DefCohSys}
Following \cite[Définition 1.1]{He}, a \textit{coherent system} of type $(n,d,k)$ on $C$ is a triple $\mathcal{E}=(E,V,\varphi)$ consisting of a coherent sheaf $E$ of rank $n$ and degree $d$, $k$-dimensional vector space $V$, and a linear map $\varphi :V \to H^0(E)$.

In addition, a morphism between coherent systems $(E,V,\varphi)$ and $(E',V',\varphi')$ is given by a linear map $f:V \to V'$ and a sheaf map $h:E \to E'$ such that the following diagram commutes 
\begin{equation} \label{equcom} \begin{split}
\xymatrix{
V \ar[r]^{\varphi} \ar[d]^f & H^0(E) \ar[d]^{H^0(h)}\\
V' \ar[r]^{\varphi'} & H^0(E')}
\end{split} \end{equation}
\end{Definition}

We also set $\rk(\eee):=n$, $\deg(\eee)=d$ and $\dim(\eee)=k$. A coherent system $\mathcal{E}=(E,V,\varphi)$ on $C$ is called \emph{pure}, if $E$ is torsion free and $\varphi$ is injective, cf. \cite[Définition 1.12]{He}. A coherent system $\cale$ is said to be \textit{injective} if $\varphi$ is injective, and we just write $\mathcal{E}=(E,V)$ understanding that $\varphi:V\hookrightarrow H^0(E)$ is the inclusion map. When $\varphi$ is an isomorphism, we say that $\eee$ is \textit{complete}, and we denote it by $|E|$. For an injective system $\eee$, we set $|\eee|:=|E|$, which stands for the complete system that has $\eee$ as a subsystem. In addition, $\cale$ is a \textit{torsion system} if $n=0$ (i.e. $E$ is either trivial or a torsion sheaf); note that every system admits a torsion subsystem whose quotient is pure. Finally, by a \textit{linear system}, we mean an injective coherent system of rank 1.

We denote the category of coherent systems on $C$ by $\csys(C)$. This is an example of a \emph{comma category}, see \cite[Example 3.2]{JNO} for details. In particular, $\csys(C)$ is abelian \cite[Main Theorem 1]{JNO} and noetherian \cite[Main Theorem 2]{JNO}; the former claim is directly stated, without proof, in \cite[p.~547 bot]{He}, while the latter is contained in \cite[Proposition 2.3]{He} as well. Furthermore, $K_0(\csys(C))=K_0(\coh(C))\oplus\Z$ by \cite[Main Theorem 1]{JNO}.

A short exact sequence in $\csys(C)$:
$$ 0 \longrightarrow \calf \longrightarrow \cale \longrightarrow \calg \longrightarrow 0. $$
where $\calf=(F,U,\phi)$, $\cale=(E,V,\varphi)$ and $\calg=(G,W,\psi)$, induces a the short exact sequence of sheaves
$$ 0 \longrightarrow F \longrightarrow E \longrightarrow G \longrightarrow 0, $$
and the following commutative diagram
\begin{equation}\label{eq:exact-diagram} \begin{split}
\xymatrix{
0 \ar[r] & U \ar[r] \ar[d]^{\phi} & V \ar[r] \ar[d]^{\varphi} & W \ar[r] \ar[d]^{\psi} & 0 \ar[d] \\ 
0 \ar[r] & H^0(F) \ar[r] & H^0(E) \ar[r] & H^0(G) \ar[r] & H^1(L) }
\end{split}\end{equation}
This allows us to define the group epimorphism 
$$ v:K_0(\csys(C))\to\Z^3 ~~,~~ v(\cale)=(\rk(E),\deg(E),\dim(V)) ~~ {\rm when } ~~ \cale\in\csys(C) ; $$
its image is called the \textit{type} of $\cale$. For any $\alpha \geq 0$, we introduce the stability function $Z_\alpha:K_0(\csys(C))\to\C$ given by, on the objects of $\csys(C)$,
$$ Z_{\alpha}(\mathcal{E})=-(d+\alpha k)+\sqrt{-1}\, n ~, $$
whose slope is
\[\mu_\alpha(\mathcal{E}) := 
\begin{cases}
    \displaystyle{\dfrac{d}{n} + \alpha\, \dfrac{k}{n}}\,\,\,\,\, \text{if $n\neq 0$} \\ ~~ \\
    \infty \,\,\,\,\,\,\,\,\,\,\,\,\,\,\,\,\,\,\,\, \, \text{if $n = 0$}
\end{cases}\]
Accordingly, we say that a coherent system $\mathcal{E}$ is $\mu_\alpha$-\emph{(semi)stable} if it is (semi)stable with respect to $Z_{\alpha}$. 

Here are some basic facts about $\mu_\alpha$-semistable systems that will be useful later on.

\begin{Proposition}\label{semistable}
Let $\mathcal{E}=(E,V,\varphi)$ be a coherent system of type $(n,d,k)$. Then 
\begin{itemize}
\item[(i)] if $n>0$ and $\mathcal{E}$ is $\mu_\alpha$-semistable, then $E$ is pure;
\item[(ii)] if $n=0$, then $\mathcal{E}$ is $\mu_\alpha$-semistable.
\item[(iii)] Every sub-system of a pure system is also pure.
\end{itemize}
\end{Proposition}
\begin{proof}
Assume $\mathcal{E}$ is $\mu_{\alpha}$-semistable and $n>0$. Set $T:={\rm Torsion}(E)$. Then we have that  $\mathcal{T}:=(T,0,0)\subset\mathcal{E}$. But $\infty=\mu_\alpha(\mathcal{T}) > \mu_\alpha(\mathcal{E})$. Thus $\mathcal{T}=0$, hence $T=0$, and so $E$ is torsion free. Also, $\mathcal{K}:=(0,\ker(\varphi),0) \subset \mathcal{E}$. But $\infty=\mu_\alpha(\mathcal{K}) > \mu_\alpha(\mathcal{E})$. Thus $\mathcal{K}=0$, hence $\ker(\varphi)=0$, and so $\varphi$ is injective. The proof is completed by \cite[Lemma 2.5]{KN}. This proves (i).

If $n=0$, for any $0 \neq \mathcal{E}'\subset \mathcal{E}$, we have $\mu_\alpha(\mathcal{E}')\leq\infty=\mu_\alpha(\mathcal{E})$. Thus $\mathcal{E}$ is $\mu_\alpha$-semistable and (ii) holds.

As for the last item, if $\calf=(F,U,\phi)$ is a subsystem of a pure system $\mathcal{E}=(E,V,\varphi)$, then $F$ is torsion-free because $E$ is, and $\varphi'$ is injective because $\varphi'\circ f = H^0(h)\circ f$ and all of the other maps in this identity are injective.
\end{proof}

\begin{Proposition}\label{prpdmo} Let $\eee\in\csys(C)$ be of type $(n,d,k)$ and $\mu_\alpha$-semistable for some $\alpha\geq 0$. If $n=0$ or $k>0$ then $d\geq 0$.
\end{Proposition}

\begin{proof}
Write $\eee=(E, V, \varphi)$. If $n=0$, then, as we have already seen, $d=h^0(E)\geq 0$. So assume $n>0$ and $k>0$. As $\eee$ is $\mu_\alpha$-semistable, then Proposition \ref{semistable}.(i) yields that $E$ is torsion-free and $k\leq h^0(E)$. Following \cite[Lemma 1.3]{Ne}, as $k>0$, let $E'$ be the subsheaf of $E$ generated by the sections of $V$. Plainly, $E'$ is torsion-free as well and, hence, of positive rank. 
Set $\eee':=(E',V,\varphi)$ which is of type, say, $(n',d',k)$. Consider the surjective map $\oc^{\oplus k}\twoheadrightarrow E'$. As $\oc^{\oplus k}$ is semistable, it follows that $0\leq d'/n'$, so $d'\geq 0$. Now $\eee'\subset \eee$, and $\eee$ is $\mu_\alpha$-semistable. So $(d'+\alpha k)/n'\leq (d+\alpha k)/n$, which yields $\big(d'+(n-n')\alpha k\big)/n'\leq d$. But $n\geq n'\geq 1$, $\alpha\geq 0$, $k>0$, $d'\geq 0$, so $d\geq 0$.
\end{proof}

Next, we want to relate $\mu_\alpha$-(semi)stability to (semi)stability. To do so, recall from \cite[Definition 2.4]{BGPMN} the notion of a \emph{virtual critical value} for a triple $(n,d,k)\in\nn^3$, and say they are $0=\alpha_0<\alpha_1 <\ldots<\alpha_N$ (according to \cite[Théorème 4.2]{He}, there are only finitely may critical values). Then one can check that \cite[Proposition 2.5]{BGPMN} applies here as well, that is, we have the following result.

\begin{Proposition}\label{stablealpha}
Let $\mathcal{E}=(E,V,\varphi)\in\csys(C)$ be of type $(n,d,k)$ and $0<\alpha <\alpha_1$. 
\begin{itemize}
\item [(i)] If $E$ is stable, then $\mathcal{E}$ is $\mu_{\alpha}$-stable;
\item [(ii)] if $\mathcal{E}$ is $\mu_\alpha$-stable, then $E$ is semistable.
\end{itemize}
\end{Proposition}

Next, we establish the existence of the Harder--Narasimhan filtration for coherent systems.

\begin{Proposition}\label{HNF}
Every nonzero coherent system on $C$ admits a Harder--Narasimhan filtration with respect to the stability function $Z_\alpha$. 
\end{Proposition}
\begin{proof}
Let $\eee\in\csys(C)$ be of rank $n$. Then $\Im(Z_\alpha(\eee))=n$. Thus, clearly, the image of $\Im(Z_\alpha)$ is $\mathbb{Z}$ which is discrete in $\mathbb{R}$. Also, by \cite[Main Theorem 2]{JNO}, $\csys(C)$ is Noetherian. So apply Proposition \ref{CHNF}.
\end{proof}

\begin{Remark} \label{rem:HNfilts}
When $\cale$ is not $\mu_\alpha$-semistable, its Harder--Narasimhan filtration will be wri-\\tten as follows:
$$ 0= \eee_0\subset \eee_1 \subset \ldots \subset \eee_l = \cale . $$
Its semistable factors are denoted by $\eee^i:=\eee_i/\eee_{i-1}$ for $i\in\{1,\ldots,l\}$; let $\cale^i=(E^i,v^i,\varphi^i)$. Also, 
$$ \mu_{\alpha, \rm max}(\eee):=\mu_\alpha(\eee^1) > \mu_\alpha(\eee^2) > \ldots > \mu_\alpha(\eee^l):=\mu_{\alpha, \rm min}(\eee) . $$
Note that we have exact sequences, for each $i=1,\dots,l$,
\begin{equation}\label{HNES}\begin{aligned}
0 &\lra \eee_{i-1} \lra \eee_{i} \lra \eee^{i} \lra 0 ,~~ {\rm and} \\
0 &\lra E_{i-1} \lra E_{i} \lra E^{i} \lra 0 .
\end{aligned} \end{equation}
Denoting the type of $\cale^i$ by $(n^i,d^i,k^i)$, we have that
$$ \sum_{i=1}^l (n^i,d^i,k^i) = (n,d,k) ~~{\rm and}~~  h^0(E) \le \sum_{i=1}^l h^0(E^i). $$
Moreover, if $\cale$ is pure, then $\cale_1=\cale^1$ is also pure, thus $\mu_{\alpha, \rm max}(\cale)<+\infty$ and each factor $\cale^i$ is also pure with $n^i>0$. In general, if $\cale$ is not pure, then $\cale_1=\cale^1$ is the maximal torsion sub-system of $\cale$ and $\cale^i$ is pure for $i\ge2$.
\end{Remark}

We conclude this section with an important result that establishes two inequalities among the parameters of a coherent system, as well as the characterization of the extremal cases. They will be essential for deriving new stability conditions on $D^b(\csys(C))$, which is the aim of Section \ref{sec:stab-tilt} below.

\begin{Theorem} \label{cota}
Let $\eee\in \csys(C)$ be a pure coherent system of type $(n,d,k)$. Then, if $\mu_{\alpha, \rm min}(\eee)\ge0$ for some $\alpha$, then $k\leq  d+n$.
\end{Theorem}

\begin{proof}
Write $\eee=(E,V,\varphi)$, and note that $\mu_\alpha(\cale)\ge\mu_{\alpha, \rm min}(\eee)\ge0$; moreover, $n>0$ since $\cale$ is pure, so $d+\alpha k\ge0$ and $d\ge0$ when $k=0$. 

Assume first that $\cale$ is $\mu_\alpha$-semistable. The proof is by induction on $n$. If $n=1$ and $k=0$, then the result follows as $d\geq 0$ if $k=0$. If $n=1$ and $k>0$, first apply Proposition \ref{prpdmo}. Then $d\geq 0$. Also, $\eee$ is $\mu_\alpha$-semistable. Therefore, by Proposition \ref{semistable}.(i), $E$ is torsion free, and of rank $1$. Thus $E$ is semistable. Again by Proposition \ref{semistable}.(i), $\varphi$ is injective. Thus $k\leq h^0(E)$. So apply Theorem \ref{thmclf1}.

Assume $n \geq 2$ and the result holds for smaller rank. If $k=0$, then the result trivially follows as seen above. So assume $k>0$. If $E$ is semistable, use the argument of the prior paragraph to conclude that $k\leq h^0(E)\le d+n$. If $E$ is not semistable, then Proposition \ref{stablealpha}.(ii) implies that $\eee$ is not $\mu_{\alpha}$-stable for $\alpha<\alpha_1$. But $\eee$ is $\mu_\alpha$-semistable for some $\alpha>0$. It follows that there exists a critical value $\alpha_*$ such that $\eee$ is $\mu_\alpha$-unstable for $\alpha<\alpha_*$ and $\mu_\alpha$-stable for $\alpha>\alpha_*$. Also, by Proposition \ref{semistable}.(i), $\varphi$ is injective. 
Then \cite[Lemma 6.5]{BGPMN} implies that $\cale$ fits into an exact sequence
\begin{equation}\label{wallcota}
    0 \to \cale_1 \to \cale \to \cale_2 \to 0
\end{equation}
where $\cale_1$ and $\cale_2$ are 
are $\mu_{\alpha_*}$-semistable coherent systems having the same $\mu_{\alpha_*}$-slope. Say $\eee_i$ is of type $(n_i,d_i,k_i)$ for $i=1,2$. First, we claim that $0<n_i<n$ for $i=1,2$. Indeed, if one of the $\eee_i$ (and hence both) has infinite $\mu_{\alpha_i}$-slope, then so does $\eee$, which implies $n=0$, a contradiction. Now assume that one of the $k_i$ vanishes, say, for instance, $k_1=0$. Then $k_2>0$ since $0<k\leq k_1+k_2$. So Proposition \ref{prpdmo} implies $d_2\geq 0$. Therefore $d_1=\big(n_1(d_2+\alpha_ik_2)\big)/n_2\geq 0$. 
So we are in a position to use induction, and hence
\[ k = k_1+k_2 \leq d_1+n_1+d_2+n_2= d+n . \]

Now assume that $\eee$ is not $\mu_\alpha$-semistable; using the notation of Remark \ref{rem:HNfilts}, note that, $n^i>0$ for all $i=1,\dots,l$ because $\cale$ is pure. Moreover, $d^i\geq 0$ when $k^i=0$ since $\mu_{\alpha, \rm min}(\eee)\geq 0$. The first part of the proof tells us that 
$k^i\leq d^i+n^i$ for every $i$. It follows that 
$$ k = \sum_{i=1}^l k^i \leq \sum_{i=1}^l d^i+n^i = d+n, $$
as desired.
\end{proof}

The next result is known in the literature as a generalization of Clifford's Theorem for vector bundles to coherent systems, compare with \cite[Theorem 2.1]{LN}.

\begin{Theorem} \label{thm:clifford-syst}
Let $\eee\in \csys(C)$ be a pure coherent system of type $(n,d,k)$.
If \linebreak $0 \leq \mu_{\alpha, \rm min}(\cale) \leq \mu_{\alpha, \rm max}(\cale)\leq 2g$ for some $\alpha$, then $k\leq  d/2+n$.
\end{Theorem}
\begin{proof}
When $\cale=(E,V,\varphi)$ is $\mu_\alpha$-semistable for some $\alpha$ and $0\leq\mu(E)\leq2g$, then \cite[Theorem 2.1]{LN} guarantees that $k\leq d/2+n$; the cited result is stated only for nonsingular curves, but the argument works for integral curves as well.

Now assume that $\cale$ is not $\mu_\alpha$-semistable and consider its Harder--Narasimhan filtration in the notation of Remark \ref{rem:HNfilts}. Since $\mu_{\alpha,\rm max}(\cale)$ is finite, we know that each $\cale^i$ is pure and each $n^i$ is positive.

We claim  $0 \leq \mu(E^i) \leq 2g$ for any $i$. Indeed, note that, for each $i=1,\dots,l$,
$$ 2g \ge \mu_\alpha(\cale^i) = \mu(E^i) + \alpha \dfrac{k^i}{n^i} \ge \mu(E^i). $$
Now, if $k^i=0$ for some $i$, then $d^i\geq 0$ since that $0\leq \mu_{\alpha, \rm min}(\cale)$ and hence $\mu(E^i)\geq 0$. If $k^i>0$, then $d^i \geq 0$ by Proposition \ref{prpdmo} and hence $\mu(E^i)\geq 0$.
Therefore, $0 \leq \mu(E^i)$, and we can conclude that $k^i\le d^i/2+n^i$.
It follows that
$$ k = \sum_{i=1}^l k^i \leq \sum_{i=1}^l \dfrac{d^i}{2}+n^i = \dfrac{d}{2}+n, $$
as desired.
\end{proof}

It will also be useful for us to characterize those $\mu_\alpha$-semistable coherent systems that satisfy the equalities in Theorem \ref{cota}.

\begin{Proposition} \label{prop:extremal}
A $\mu_\alpha$-semistable coherent system $\eee$ of type $(n,d,k)$ with $n,k>0$ sa-\\tisfies $k=d+n$, if and only if either $\eee=|\oc^{\oplus n}|$ when $g\ge1$ or $\eee=\oplus_{i=1}^n\left| \mathcal{O}_{\mathbb{P}^1}(a_i) \right|$, $a_i\geq 0$, when $g=0$.
\end{Proposition}
\begin{proof}
Assume first that $\eee$ satisfies $k=d+n$; Proposition \ref{prpdmo} also guarantees that $d\ge0$.

If $\mu(E)\le2g$, then $d+n\le d/2+n$  by \cite[Theorem 2.1]{LN}, thus $d\le0$; it follows that $d=0$, and $k=n$. We claim that $\eee=|\oc^{\oplus n}|$. If $E$ is semistable, the result follows from Proposition \ref{charsheaf}. Now,  If $E$ is not semistable, then Proposition \ref{stablealpha}.(ii) implies that $\eee$ is not $\mu_{\alpha}$-stable for $\alpha<\alpha_1$. But $\eee$ is $\mu_\alpha$-semistable for some $\alpha>0$. Hence, we follow the proof of Theorem \ref{cota}.  From \cite[Lemma 6.1]{BGPMN} we have that in the exact sequence (\ref{wallcota}), $k_1/n_1<k/n=1$ which implies $k_1<n_1$. Since $\eee_i$ is $\mu_\alpha$-semistable for $i=1,2$, it follows that $k=k_1+k_2 <n_1+n_2=n$ which is a contradiction.

When $\mu(E)\ge2g$, then \cite[Theorem 2.1]{LN} implies that $d+n\leq d+n(1-g)$, so in fact $g=0$. In this case $\eee=\oplus_{i=1}^n\left| \mathcal{O}_{\mathbb{P}^1}(a_i) \right|$, $a_i\geq 0$ by \cite[Lemma 3.1]{LNg0}
\end{proof}

\begin{Remark}
\emph{Note that the hypotheses in the above result are necessary. In fact, for $n=0$, whatever is $C$, consider $\eee=(0,V,0)$ with $\dim(V)=k>0$. Then $n$ and $d$ vanish. Also, $\eee$ is $\mu_\alpha$-semistable for any $\alpha$ since $\mu_\alpha(\eee)=\infty$. So $k>0=d+n$ contradicts the statement of Lemma \ref{cota}. As for $n>0$, $d<0$ and $k=0$, consider $\eee=(\oo_{\mathbb{P}^1}(-a),0,0)$ with $a\geq 2$ on $C=\mathbb{P}^1$. Then $n=1$, $d=-a$ and $k=0$.  Also, $\eee$ is $\mu_\alpha$-semistable for any $\alpha$ since if $\eee'\subset \eee$ is of type $(n',d',k')$ then $n'=1$, $d'\leq -a$ and $k'=0$, hence $\mu_\alpha(\eee')\leq\mu_\alpha(\eee)$. Thus $k=0>-a+1=d+n$ gainsays Lemma \ref{cota} as well.} 
\end{Remark}


\section{Standard Stability Conditions for Coherent Systems} \label{sec:stab-cs}

We already noticed that $(\csys(C),Z_\alpha)$ is a stability condition in the sense of Definition \ref{dfnstc}, since $Z_\alpha$ is a stability function and Proposition \ref{HNF} provides a Harder--Narasimhan filtration for each $\cale\in\csys(C)$. However, for later use, we will introduce another fami-\\ly of stability functions on $\csys(C)$, depending on a new parameter, and of which the one given above is a particular case. Namely, for $\eee\in\csys(C)$ of type $(n,d,k)$, set
\begin{equation}
\label{equzab}
Z_\alpha^\beta(\eee):=Z_\alpha(\eee)+n\beta=-(d+\alpha k-\beta n) +\sqrt{-1}\, n
\end{equation}
So $Z_{\alpha}$ corresponds to the case $\beta=0$. Clearly, the slope with respect to $Z_\alpha^\beta$ is given by
\begin{equation}
\label{equmim}
\mu_\alpha^\beta(\eee)=\mu_\alpha(\eee)-\beta
\end{equation}
Therefore, a system is $\mu_\alpha^\beta$-(semi)stable if and only if $\mu_\alpha$-(semi)stable, so the Harder--Narasimhan filtration with respect to $Z_\alpha$ is also the Harder--Narasimhan filtration with respect to $Z_\alpha^\beta$.

Since $\csys(C)$ can also be regarded as the heart of the standard t-structure on the derived category $D^b(\csys(C))$, the pair $(\csys(C),Z_\alpha^\beta)$ can also be regarded as a two-parameter family of pre-stability conditions on $D^b(\csys(C))$. We will now check that they satisfy the support property.

\begin{Proposition}
\label{prpsupSC}
For any $\alpha>0$, the following holds: 
\begin{equation} \label{equsup}
S_\alpha:=\inf\left\{\frac{|Z_\alpha (\mathcal{E})|}{\left \| v(\mathcal{E}) \right \|}\ \bigg{|}\ 0 \neq \mathcal{E} \in \csys(C) \,\,  \text{is $\mu_\alpha$-semistable}\right\} > 0,
\end{equation}
that is, the stability function $Z_\alpha$ satisfies the support property.
\end{Proposition}

\begin{proof}
Given $\mathcal{E} \in \csys(C)$ $\mu_\alpha$-semistable of type $(n,d,k)$, we have
$$ \frac{|Z_\alpha (\mathcal{E})|^2}{\left \| v(\mathcal{E}) \right \|^2}= 
\frac{(d+\alpha k)^2+ n^2}{n^2+d^2+k^2}=:f(n,d,k) $$ 
Now if $\alpha \leq 1$, then $f(n,d,k)\geq \alpha^2$. Indeed, the inequality reads  $2d\alpha k\geq (\alpha^2-1)(n^2+d^2)$, which holds as $0<\alpha\leq 1$, and as either $k= 0$ or if $k>0$, then $d\geq 0$ by Proposition \ref{prpdmo} since $\eee$ is $\mu_\alpha$-semistable. And if $\alpha> 1$ then $f(n,d,k)\geq 1$. In this case, the inequality reads $2d\alpha k \geq (1-\alpha^2)k^2$, which holds for similar reasons. Thus \eqref{equsup} holds.
\end{proof}

Following Proposition \ref{prposp}, the symmetric bilinear form $Q_\alpha$ on $\Lambda_\R=\R^3$ associated to the support property above is such that
\begin{equation}
\label{eququa}
Q_\alpha(n,d,k) = 
\begin{cases}
    \displaystyle{\frac{(1-\alpha^2)(n^2+d^2)+2\alpha dk}{\alpha^2}}, \,\,\,\,\, \text{if $0<\alpha\leq 1$}, \\ ~~ \\
  (\alpha^2-1)k^2+2\alpha dk, \,\,\,\,\,\,\,\,\,\,\,\,\,\,\,\,\,\,\, \, \ \  \text{if $\alpha\geq 1$.}
\end{cases}
\end{equation}

We argue that the same quadratic form satisfies the properties of Proposition \ref{prposp}.

\begin{Lemma} \label{Quadratic}
Consider the symmetric bilinear form $Q_\alpha$ on $\Lambda_\mathbb{R}$. Then
\begin{itemize}
\item [(i)] $Q_\alpha(\eee)\geq 0$ for all $\mu_{\alpha}^\beta$-semistable $\mathcal{E} \in \csys(C)$;
\item [(ii)] $Q_\alpha(w)< 0$ for all nonzero $w\in \Lambda_{\mathbb{R}}$ with $Z_{\alpha}^\beta(w)=0$.
\end{itemize}
\end{Lemma}
\begin{proof}
To prove (i), note that, by \eqref{equmim}, $\eee$ is $\mu_\alpha^\beta$-semistable if and only if it is $\mu_\alpha$-semistable. So the claim follows from the proof of Proposition \ref{prpsupSC}.  For (ii), let $w=(n,d,k)\in\Lambda_\mathbb{R}$ be a nonzero vector. If $Z_{\alpha}^\beta(w)=0$, then $n=0$ and $d=-\alpha k$. Therefore \eqref{eququa} yields $Q_\alpha(w)=-(\alpha^2+1)k^2$. But $k\neq 0$ since $w\neq 0$. So $Q_\alpha(w)<0$.
\end{proof}

By putting all the pieces together, we have proven the following result.

\begin{Theorem} \label{thmzab} 
$(\csys(C), Z_\alpha^\beta)$ is a stability condition on $D^b(\csys(C))$ for each
$\alpha\geq0$ and $\beta\in\R$. Moreover, the map 
\[ \begin{matrix}
\mathbb{R}_{\geq 0}\times \mathbb{R} & \to & 
\stab(D^b(\csys(C)))\\
(\alpha,\beta) & \longmapsto & (\csys(C), Z_\alpha^\beta)
\end{matrix} \]
is a continuous embedding.
\end{Theorem}

The stability conditions described in Theorem \ref{thmzab} are called \textit{standard stability conditions}.

We conclude this section by studying the orbits of $(\csys(C), Z_\alpha^\beta)$ under the $\widetilde{GL}^+(2,\mathbb{R})$ action described in Section \ref{sec:stab-tri}.

\begin{Lemma} \label{action1}
Consider two pairs $(\alpha,\beta)$ and $(\alpha',\beta')$ in $\mathbb{R}_{\geq 0}\times \mathbb{R}$.  If $\alpha\ne\alpha'$, then there is no $T\in GL^+(2,\R)$ such that $Z_\alpha^\beta = T\circ Z_{\alpha'}^{\beta'}$.
\end{Lemma}

Therefore, the image of the embedding provided in Theorem \ref{thmzab} is not contained in a single $\widetilde{GL}^+(2,\mathbb{R})$-orbit.

\begin{proof}
Let $T= \begin{pmatrix} A & B \\ C & D \end{pmatrix} \in GL^+(2,\mathbb{R})$. If $Z_\alpha^\beta = T\circ Z_{\alpha'}^{\beta'}$, then
$$ A=1 ~~,~~ A\alpha = \alpha' ~~, -A\beta+B=\beta' ~~, $$
$$ C=0 ~~, D=1. $$
In particular, this system of equations does not admit any solutions when $\alpha\ne\alpha'$.
\end{proof}

Note, however, that
$\left(\begin{array}{cc} 1 & \beta-\beta' \\ 0 & 1 \end{array} \right) Z_\alpha^\beta = Z_\alpha^{\beta'}$.


\section{Stability Conditions for Tilted Coherent Systems} \label{sec:stab-tilt}

In this section, we construct a three-parameter family of stability conditions on the $D^b(\csys(C))$. To achieve this, we use the tilting construction described in Section \ref{sec:tilt} applied to the category of coherent systems and the stability function $Z_\alpha$ given in the previous section. 

Given $\alpha,\beta\in\mathbb{R}$, with $\alpha \geq 0$, set, following the \eqref{equtbt} and \eqref{equfbt}:
\[ \begin{aligned}
\tab := \mathscr{T}^\beta(\csys(C),Z_{\alpha})&=\{\mathcal{E}\in \csys(C)\, |\, \text{$\mu_\alpha(\mathcal{G}) > \beta$ whenever $\mathcal{E} \twoheadrightarrow \mathcal{G}$}\} \\
\fab := \mathscr{F}^\beta(\csys(C),Z_{\alpha}) & =\{\mathcal{E}\in \csys(C)\, |\, \text{ $\mu_\alpha(\mathcal{F}) \leq \beta$  for every $0\neq\mathcal{F} \subset\mathcal{E}$}\}
\end{aligned} \]
Propositions \ref{Homomorphism} and \ref{HNF} show that $(\mathcal{T}_\alpha^\beta, \mathcal{F}_\alpha^\beta)$ is a torsion pair on $\csys(C)$. Set
$$ \csys^\beta_\alpha(C) = \langle\mathcal{F}^\beta_\alpha[1], \mathcal{T}^\beta_\alpha\rangle ; $$
in other words, $\csys_\alpha^\beta(C)$ is the tilt of $\csys(C)$ with respect to the torsion pair $(\mathcal{T}_{\alpha}^\beta,\mathcal{F}_\alpha^\beta)$. Proposition \ref{newheart} yields that $\csys^\beta_\alpha(C)$ is the heart of a bounded $t$-structure on $D^b(\csys(C))$. In addition, one can check that 
\[ \csys^\beta_\alpha(C) := \bigg{\{} \mathcal{E} \in D^b(\csys(C)) \, \bigg{|}\, \begin{matrix}\ \hh^i(\mathcal{E}) =0 \,\,\, \text{for all $i \neq 0,-1$}, \\ \ \hh^0(\mathcal{E}) \in \mathcal{T}^\beta_\alpha\ \text{and}\ \hh^{-1}(\mathcal{E}) \in \mathcal{F}^\beta_\alpha \end{matrix} \bigg{\}}. \]

Now we will define a stability condition on $\csys^\beta_\alpha(C)$. We will start by setting a map on $\csys(C)$, which we will extend in \eqref{equTOM} to get a central charge in the tilted category.


Given $\alpha,\beta,\gamma\in\mathbb{R}$ with $\alpha\ge0$ and $\mathcal{E}\in \csys(C)$ of type $(n,d,k)$, set
\begin{equation}
\label{equctc}
Z^{\beta, \gamma}_\alpha(\mathcal{E})= (d+\gamma n-k) + \sqrt{-1}\,(d+ \alpha \, k - \beta \, n).
\end{equation}
Now let $\mathcal{E}\in \csys^\alpha_{\beta}(C)$. Then $\mathcal{E}$ fits into an exact triangle
\begin{equation}
\label{equTRI}
0 \to \hh^{-1}(\mathcal{E})[1] \to \mathcal{E} \to \hh^{0}(\mathcal{E}) \to 0
\end{equation}
where $\hh^{-1}(\mathcal{E}) \in \mathcal{F}^\beta_\alpha$ and $\hh^{0}(\mathcal{E})\in \mathcal{T}^\beta_\alpha$. Set
\begin{equation}
\label{equTOM}
Z^{\beta, \gamma}_\alpha(\mathcal{E}) := Z^{\beta, \gamma}_\alpha(\hhz(\mathcal{E}))-Z^{\beta, \gamma}_\alpha(\hhm(\mathcal{E}))
\end{equation}
and extend this map additively to get a homomorphism of abelian groups
$$ Z_{\alpha}^{\beta,\gamma}:K_0(\csys_\alpha^\beta(C))\to \mathbb{C}, $$
which factors through the lattice $\Z^3$.

\begin{Proposition}\label{slopes}
Let $\mathcal{E}\in \csys(C)$ be $\mu_\alpha$-semistable of type $(n,d,k)$. Assume that $d\geq 0$ if $k=0$ and $d\geq k$ if $n=0$. Then the following hold:
\begin{itemize}
\item[(i)] if $n=0$, then $\mu^{\beta, \gamma}_\alpha(\mathcal{E})\leq 0$, with equality if and only if $d=k$;
\item[(ii)] if $n>0$ and $\mu_{\alpha}(\eee)>\beta$, then $\mu^{\beta, \gamma}_\alpha(\mathcal{E})<0$;
\item[(iii)] if $n>0$ and $\mu_{\alpha}(\eee)<\beta$, then $\mu^{\beta, \gamma}_\alpha(\mathcal{E})>0$.
\end{itemize}
\end{Proposition}

\begin{proof}
If $n=0$, then $\mu^{\beta, \gamma}_\alpha(\mathcal{E})=-(d-k)/(d+\alpha k)$ and (i) follows since $d\geq k$ or $k=0$. 

If $n>0$, write
\begin{equation}
\label{equmn>}
{\mu^{\beta, \gamma}_\alpha(\mathcal{E})=-\,\frac{\displaystyle{\frac{d-k}{n}+\gamma}}{\mu_{\alpha}(\eee)-\beta}}
\end{equation}
and (ii) and (iii) follow from Lemma \ref{cota}.
\end{proof}

\begin{Proposition}\label{TorHNF}
Let $\mathcal{E} =(E, V, \varphi)\in \csys(C)$ be of type $(n,d,k)$. The following holds:
\begin{itemize}
   \item [(i)] if $\mathcal{E} \in \mathcal{F}_\alpha^\beta$, then $n>0$ and $\eee$ is pure and injective;
  \item [(ii)] if $\mathcal{E} \in \mathcal{T}_\alpha^\beta$, then
 \begin{itemize}
    \item [(a)] if $\eee$ is pure, then $d \geq 0$, $k \leq d+n$;
   \item [(b)] if $\eee$ is not pure $n>0$, then there is $\mathcal{T}\subset\eee$ of type $(0,t,u)$ such that $\eee/\mathcal{T}$ is pure and $k \leq d-t+n+u$; in particular, if $\eee$ is injective, then $k \leq d+n$.
\end{itemize} 
\end{itemize}
\end{Proposition}

\begin{proof}
   Write $\eee=(E,V,\varphi)$ and set $T:={\rm Torsion}(E)$. To prove (i), if $n=0$, then we have $\beta <\infty=\mu_\alpha(\mathcal{E})$ which is precluded since $\mua(\eee)\leq\beta$ as $\mathcal{E}\in \mathcal{F}_\alpha^\beta$. Let $\mathcal{T}:=(T,0,0)$ and assume $T\neq 0$. Then $\rank(T)=0$ and $\mua(\mathcal{T})=\infty$. But $\mathcal{T} \subset \eee\in\fab$, so $\mua(\mathcal{T})\leq \beta$, which is a contradiction. Thus $T=0$, hence $E$ is torsion free, and $\mathcal{E}$ is pure.  Finally, set $\mathcal{K}:=(0, \ker(\varphi),0)$ and assume $\ker(\varphi)\neq 0$. Then $\mua(\mathcal{K})=\infty$. But $\mathcal{K} \subset \eee$, so $\mua(\mathcal{T})\leq \beta$, which is again a contradiction. Thus $\ker(\varphi)=0$, and $\mathcal{E}$ is injective.
   
To prove (ii).(a), let $\eee^i$ be the factors of $\eee$. From Remark \ref{filtrations}, $\eee^i \in \tab$ for all $i$. Also, $\eee^i$ is $\mu_\alpha$-semistable. Say $\eee^i$ is of type $(n^i,d^i,k^i)$. If $k_i=0$, as $\eee^i \in \tab$, we have $\mua(\eee^i)=d^i/n^i>\beta>0$, so $d^i>0$. If $k^i>0$, then $d^i\geq 0$ by Proposition \ref{prpdmo}. So $d=\sum d_i\geq 0$. Also, $\mu_{\alpha, \rm min}(\eee)=\mua(\eee^k)>\beta>0$. So $k\leq d+n$ by Theorem \ref{cota}.

To prove (ii).(b), write $\eee=(E,V,\varphi)$ and view $H^0(T)\subset H^0(E)$. Set $U:=\varphi^{-1}(H^0(T))$. Consider the system $\mathcal{T}=(T,U,\varphi|_U)$ and say $\mathcal{T}$ is of type $(0,t,u)$. So $\eee/\mathcal{T}$ is of type $(n,d-t,k-u)$. As $\eee/\mathcal{T} \in \tab$ and is pure, Then $k\leq d-t+n-u$ by (ii)(a). If $\eee$ is injective, then $\mathcal{T}$ is injective which implies $t\leq u$, and hence $k\leq d+n-(t-u)\leq d+n$. 
\end{proof}

\begin{Proposition} \label{WSC}
If $\beta\geq0$ and $\gamma>1$, then $Z^{\beta,\gamma}_\alpha: K_0(\csys_\alpha^\beta(C))\to\mathbb{C}$ is a stability function on $\csys_\alpha^\beta(C)$. 
\end{Proposition} 
\begin{proof} Let $\mathcal{E}\in \csys^\alpha_{\beta}(C)$. As we have just seen, $\mathcal{E}$ fits into an exact 
    \[ 0 \to \hh^{-1}(\mathcal{E})[1] \to \mathcal{E} \to \hh^{0}(\mathcal{E}) \to 0\]
    where $\hh^{-1}(\mathcal{E}) \in \mathcal{F}^\beta_\alpha$ and $\hh^{0}(\mathcal{E})\in \mathcal{T}^\beta_\alpha$. Say $\hh^{-1}(\mathcal{E})$ is of type $(n_1,d_1,k_1)$ and $\hh^{0}(\mathcal{E})$ of type $(n_0,d_0,k_0)$. Then
    \[
    \begin{aligned}
        \Im(Z^{\beta, \gamma}_\alpha(\mathcal{E})) &= \Im(Z^{\beta, \gamma}_\alpha(\hh^{-1}(\mathcal{E})[1])) + \Im(Z^{\beta, \gamma}_\alpha(\hh^{0}(\mathcal{E}))) \\
        &= -\Im(Z^{\beta, \gamma}_\alpha(\hh^{-1}(\mathcal{E}))) + \Im(Z^{\beta, \gamma}_\alpha(\hh^{0}(\mathcal{E})))\\
        & =-(d_1+ \alpha\, k_1 - \beta\, n_1)+ (d_0+ \alpha\, k_0 - \beta\, n_0)
    \end{aligned}\]
As $\hh^{-1}(\mathcal{E}) \in \mathcal{F}^\beta_\alpha$, then $n_1>0$ by Proposition \ref{TorHNF}(i). Assume $n_0>0$ too. Then
$$
\Im(Z^{\beta, \gamma}_\alpha(\mathcal{E}))=n_1(\beta-\mu_\alpha(\hh^{-1}(\mathcal{E})))+n_0(\mu_\alpha(\hh^{0}(\mathcal{E}))-\beta)
$$
and hence $\Im(Z^{\beta, \gamma}_\alpha(\mathcal{E}))>0$. So assume $n_0=0$. Then
$$
\Im(Z^{\beta, \gamma}_\alpha(\mathcal{E}))=n_1(\beta-\mu_\alpha(\hh^{-1}(\mathcal{E})))+d_0+ \alpha\, k_0
$$
which remains non-negative and vanishes if and only if $\mu_\alpha(\hh^{-1}(\mathcal{E}))=\beta$ and $d_0=k_0=0$. As $n_0=0$,  we have $\hh^{0}(\mathcal{E})=0$. 

We therefore conclude that $\Im(Z^{\beta, \gamma}_\alpha(\mathcal{E}))\ge0$, and equality holds if and only if $\hh^{0}(\mathcal{E})=0$ and $\mu_\alpha(\hh^{-1}(\mathcal{E}))=\beta$, which implies that $\hh^{-1}(\mathcal{E})$ is $\mu_\alpha$-semistable by Lemma \ref{prpelm}(iii). Thus $\Im(Z^{\beta, \gamma}_\alpha(\mathcal{E}))=0$ implies that
\begin{equation} \label{equrez}   
\Re(Z^{\beta, \gamma}_\alpha(\mathcal{E}))=
-\Re(Z^{\beta, \gamma}_\alpha( \hh^{-1}(\mathcal{E}))) = -(d_1+\gamma n_1 -k_1) < -(d_1+n_1-k_1) \le 0.  
\end{equation}
where the first inequality holds because $\gamma>1$ and the second holds due to  Theorem \ref{cota} provided $\beta\ge0$.
\end{proof}

We record, for later reference, the following statement that was established as part of the previous proof. 
\begin{Corollary}\label{Im=0}
      Let  $\mathcal{E}\in \csys^\alpha_{\beta}(C)$
      be such that $\Im(Z^{\beta, \gamma}_\alpha(\mathcal{E}))=0$. Then  $\hh^{0}(\mathcal{E})=0$, $\mathcal{E} \cong \hh^{-1}(\mathcal{E})[1]$ and $\hh^{-1}(\mathcal{E})$ is $\mu_\alpha$-semistable with $\mu_\alpha(\hh^{-1}(\mathcal{E}))=\beta$.  
  \end{Corollary}

The next step is to show that the stability function $Z^{\beta,\gamma}_\alpha$ defines a stability condition on $\csys_\alpha^\beta(C)$. The approach is based on Proposition \ref{CHNF}, so we need a few technical lemmas.

\begin{Lemma}\label{discrete}
The following holds.
\begin{itemize}
\item [(i)] The image of $\Im(Z^{\beta, \gamma}_\alpha)$ is discrete in $\mathbb{R}$ if and only if $\alpha, \beta \in \mathbb{Q}$. 
\item [(ii)] The image of $\Re(Z^{\beta, \gamma}_\alpha)$ is discrete in $\mathbb{R}$ if and only if $\gamma \in \mathbb{Q}$.
\end{itemize}
\end{Lemma}

\begin{proof}
Set $A:=\{d+\alpha k-\beta n\, |\, (n,d,k) \in \mathbb{Z} ^3\}$, which agrees with the image of $\Im(Z^{\beta, \gamma}_\alpha)$. If $\alpha,\beta\in \mathbb{Q}$, write $\alpha=a_\alpha/b_\alpha$ and $\beta=a_\beta/b_\beta \in \mathbb{Q}$. Fix $a_0=d_0+\alpha k_0-\beta n_0 \in A$ and let $a=d+\alpha k-\beta n\in A$ be arbitrary with $a\neq a_0$. Then
\begin{align*}
|a_0-a| & = |(d_0-d)+\alpha(k
   _0-k)+\beta(n-n_0)|\\
   &=\frac{|b_\alpha b_\beta(d_0-d)+a_\alpha b_\beta(k
_0-k)-a_\beta b_\alpha(n_0+n)|}{|b_\alpha b_\beta|}\geq \frac{1}{|b_\alpha b_\beta|}
\end{align*}
where the last inequality holds as the numerator is a natural number, and also nonzero as $a\neq a_0$. Therefore $A$, that is, the image of $\Im(Z^{\beta, \gamma}_\alpha)$, is discrete.

Conversely, assume $\alpha$ is not rational (a similar argument holds for $\beta$). Then Dirichlet's Approximation Theorem asserts that for any real numbers $\alpha$ and $N\geq 1$, there are integers $p$ and $q$, with $1\leq q\leq N$, such that $|q\alpha-p|<1/N$. So given $N\in\mathbb{N}^*$, take $k=q$, $d=-p$ and $n=0$. Now vary $N$ and note $|q\alpha -p|$ never vanishes as $\alpha$ is irrational. This shows that $0$ is not an isolated point of $A$. Thus the image of $\Im(Z^{\beta, \gamma}_\alpha)$ is not discrete. This proves (i). The proof of (ii) is similar.
\end{proof}

\begin{Lemma}\label{Noether}
If $\alpha,\beta\in\mathbb{Q}$, then the tilted category $\csys^\beta_\alpha(C)$ is Noetherian.
\end{Lemma}

\begin{proof}
Let $\mathcal{E} \in \csys^\beta_\alpha(C)$ and let 
\begin{equation}
\label{equinc}
0=\mathcal{E}_0 \subset \mathcal{E}_1 \subset \cdots \subset \mathcal{E}_l \subset \cdots \subset \mathcal{E}
\end{equation}
be a chain of subobjects of $\csys^\beta_\alpha(C)$.
Consider the exact sequence 
\begin{equation}
\label{eququo}
0\to  \mathcal{E}_l \to  \mathcal{E} \to \mathcal{E}/\mathcal{E}_l \to 0.
\end{equation}
By Proposition \ref{WSC}, $\Im(Z_\alpha^{\beta, \gamma})$ is a non-negative function. Thus $\Im(Z_\alpha^{\beta, \gamma})(\mathcal{E}_l)\leq \Im(Z_\alpha^{\beta, \gamma}(\mathcal{E})$ for all $l$. Also, by Lemma \ref{discrete}, $\Im(Z_\alpha^{\beta, \gamma})$ is discrete, so we have finitely many possibilities for $\Im(Z_\alpha^{\beta, \gamma}(\mathcal{E}_l))$. On the other hand, for every $l$, we have an exact sequence
$$0\to \mathcal{E}_{l-1} \to  \mathcal{E}_l \to \mathcal{E}_l/\mathcal{E}_{l-1} \to 0. $$
Thus, again, $\Im(Z_\alpha^{\beta, \gamma})(\mathcal{E}_{l-1})\leq \Im(Z_\alpha^{\beta, \gamma})(\mathcal{E}_l)$. In other words,  $\Im(Z_\alpha^{\beta, \gamma}(\mathcal{E}_l))$ is an increasing function of $l$ and, from the above, it reaches finitely many values.  Thus, it stabilizes, say at $l_0$. So we may start the inclusions in \eqref{equinc} by $l_0$. Also, we may replace $\eee_l$ by $\eee_l/\eee_{l_0}$ for all $l\geq l_0$, as if the latter stabilizes, so does the former. Now the exact sequences
$$
0\to \mathcal{E}_{l_0} \to  \mathcal{E}_l \to \mathcal{E}_l/\mathcal{E}_{l_0} \to 0.
$$
yield that we may further assume $\Im(Z_\alpha^{\beta, \gamma}(\mathcal{E}_l))=0$ for all $l$. Therefore, as it was discussed in the proof of Proposition \ref{WSC}, $\mathcal{H}^0(\mathcal{E}_l)=0$ for all $l$. Set $Q_l:=\mathcal{E}/\mathcal{E}_l$. We have that \eqref{eququo} induces the  exact sequences 
\begin{equation} \label{equhmu}
0 \to \hh^{-1}(\mathcal{E}_l) \rightarrow \hh^{-1}(\mathcal{E}) \rightarrow \hh^{-1}(\mathcal{Q}_l) \rightarrow 0
\end{equation}
\begin{equation} \label{equhze}
 0 \to  \hh^{0}(\mathcal{E}) \rightarrow \hh^{0}(\mathcal{Q}_l) \rightarrow 0
\end{equation}
in $\csys(C)$. Thus \eqref{equhze} yields $\hh^{0}(\mathcal{Q}_l)=\hh^{0}(\mathcal{E})$ for any $l$. Now the sequence 
\[
    0 = \hh^{-1}(\mathcal{E}_0)  \hookrightarrow \hh^{-1}(\mathcal{E}_1)  \hookrightarrow \hh^{-1}(\mathcal{E}_2)  \hookrightarrow \ldots \hookrightarrow \hh^{-1}(\mathcal{E})
\]
stabilizes because $\csys(C)$ is Noetherian.  But \eqref{equhmu} yields $\hh^{-1}(\mathcal{Q}_l) = \hh^{-1}(\mathcal{E})/\hh^{-1}(\mathcal{E}_l)$. Therefore, for sufficiently large $l$, we have that  $\hh^{-1}(\mathcal{Q}_l)$ is constant as well. It follows that $\mathcal{Q}_l$ stabilizes and so does $\mathcal{E}_l$, as desired.
\end{proof}

\begin{Proposition}\label{HNFC}
If $\alpha, \beta \in \mathbb{Q}$, then every nonzero object in $\csys^\beta_\alpha(C)$ admits a Harder--Narasimhan filtration with respect to $Z^{\beta, \gamma}_\alpha$. 
\end{Proposition}

\begin{proof}
    Since $\csys^\beta_\alpha(C)$ is Noetherian by Lemma \ref{Noether}, and the image of $\Im(Z^{\beta, \gamma}_\alpha)$ is discrete in $\mathbb{R}$ by Lemma \ref{discrete}, the statement follows from Proposition \ref{CHNF}. 
 \end{proof} 

With the technicalities resolved, we are finally ready to state the main result of this section.

\begin{Theorem} \label{StabCond}
$\tau_{\alpha}^{\beta,\gamma}:=(\csysab(C),Z_\alpha^{\beta,\gamma})$ is a pre-stability condition in $D^b(\csys(C))$ if $\alpha, \beta \in \mathbb{Q}$ with $\alpha \geq0$ and $\beta\ge0$, and every $\gamma>1$. If, in addition, $\gamma\in\Q$, then $\tau_\alpha^{\beta,\gamma}$ is locally finite.
\end{Theorem}

\begin{proof}
From Proposition \ref{WSC}, $Z_\alpha^{\beta,\gamma}$ is a stability function on $\csysab(C)$, and from Proposition \ref{HNFC}, every nonzero object in $\csysab(C)$ has a Harder--Narasimhan filtration with respect to $Z_\alpha^{\beta,\gamma}$. Thus $(\csysab(C),Z_\alpha^{\beta,\gamma})$ is a pre-stability condition in $D^b(\csys(C))$ under the stated conditions on the parameters $(\alpha,\beta,\gamma)$.

If $\alpha,\beta,\gamma\in\Q$, then the image of $Z_\alpha^{\beta,\gamma}$ is discrete by Lemma \ref{discrete}. Then local finiteness follows from the proof of \cite[Lemma 4.4]{Bridgeland1}.
\end{proof}

Establishing the support property for $\tau_\alpha^{\beta,\gamma}$ will require further work, which we will undertake in the next section.


\section{Bogomolov--Gieseker Inequality for Coherent Systems} \label{sec:bg-ineq}

In this section, we will prove that the pre-stability conditions $(\csys_\alpha^\beta(C), Z_\alpha^{\beta,\gamma})$ satisfy the support property
under certain weak impositions on their parameters.

We will approach the entire problem from a different framework, as outlined in Piyaratne and Toda \cite[Section 2]{PT}, which we will adapt to our specific context. The key ingredient is finding a suitable \textit{Bogomolov-Gieseker inequality} (in the sense of \cite{PT} for $\mu_\alpha$-semistable coherent systems. This inequality will allow us to achieve the desired results for the tilted category of coherent systems.

Let $\sigma=(\cata,Z)$ be a stability condition on a triangulated category $\catd$ such that $\im(Z)\subset \mathbb{Q}+\sqrt{-1}\,\mathbb{Q}$. Following \cite[ Definition 2.7]{PT}, we say $\sigma$ satisfies a \emph{Bogomolov--Gieseker (BG) inequality} if there are linear maps
$$ \Delta_R,\Delta_I:\Lambda\otimes\mathbb{Q}\to\mathbb{Q} $$
such that for any semistable object $E\in\cata$, the following inequality holds:
\begin{equation} \label{bg-ineq}
\Delta(E) := \Re(Z(E))\cdot\Delta_R(E) + \Im(Z(E))\cdot\Delta_I(E) \ge0.
\end{equation}
Let $M$ be the slope with respect to $Z$, and consider the torsion pair 
\begin{equation}
\label{torpar}
\begin{aligned}
    \mathscr{T}(\cata,Z)& :=\{E\in \cata\, |\, \text{$M(G) > 0$ whenever $E \twoheadrightarrow G$}\} \\
    \mathscr{F}(\cata,Z)& :=\{ E\in \cata\, |\, \text{ $M(F) \leq 0$  for every $0\neq F \subset E$}\}
    \end{aligned}
\end{equation}
and the tilted category
\begin{equation} \label{eq:tilt}
\cata^\dagger := \langle \mathscr{F}[1],\mathscr{T} \rangle.
\end{equation}
In \cite[p.~9]{PT} introduces a family of homomorphisms $Z^\dagger_t:\Lambda\to\mathbb{C}$, for $t\ge0$, given by
\begin{equation} \label{eq:zdagger}
Z^\dagger_t := -\sqrt{-1}\cdot Z +t\cdot \Delta_I = (\Im(Z)+t\Delta_I) -\sqrt{-1}\cdot \Re(Z)
\end{equation}

In \cite[Lemma 2.11]{PT} it is proved that $Z^\dagger_t$ is a stability function on $\cata^\dagger$. Moreover, $\cata^\dagger$ is Noetherian by \cite[Lemma 2.17]{PT}, and the Harder--Narasimhan property holds for $(\cata^\dagger,Z^\dagger_t)$ by \cite[Lemma 2.18]{PT}.
Finally, \cite[Corollary 2.22]{PT} assures that if $\sigma=(\cata,Z)$ satisfies the BG inequality \eqref{bg-ineq}, then $\sigma^\dagger=(\cata^\dagger,Z^\dagger_t)$, $t>0$, is a stability condition on $\catd$; the support property is verified by the quadratic form 
\begin{equation} \label{eq:qdagger}
Q^\dagger_{r,t} := r\cdot Q + t\cdot\Delta
\end{equation}
where $Q$ is the quadratic form arising from the support property of $\sigma=(\cata,Z)$, $0<r<r_0$, where $r_0$ is determined by $Q$, and $t>0$ \cite[Lemma 2.20]{PT}. Moreover, the map
$$ \mathbb{R}_{>0} \to \stab(\catd) ~~{\rm given~by} ~~
t \mapsto (\cata^\dagger,Z_t^\dagger) $$
is continuous \cite[Corollary 2.22]{PT}.

Let us now adjust this construction to our context. For $\alpha\in\mathbb{Q}_{> 0}$ and $\beta\in\mathbb{Q}_{\ge0}$, given $\eee\in\csys(C)$ of type $(n,d,k)$, recall \eqref{equzab} and set
\begin{equation}
\label{equnab}
Z_\alpha^\beta(\eee) := Z_\alpha(\eee) + \beta n = -(d + \alpha k-\beta n) +\sqrt{-1}\cdot n.
\end{equation}
As we have already seen, it yields a stability function $Z_\alpha^\beta:K_0(\csys(C))\to\mathbb{C}$ factoring through the lattice $\Lambda=\mathbb{Z}^3$. Also, $\im(Z_\alpha^\beta)\subset\mathbb{Q}+\sqrt{-1}\,\mathbb{Q}$ as $\alpha,\beta\in\mathbb{Q}$. By Theorem \ref{thmzab}, $(\csys(C),Z_\alpha^\beta)$ is a stability condition on $D^b(\csysab)$.

Let $\mu_\alpha^\beta$ be the slope with respect to $Z_\alpha^\beta$. Then $\mu_\alpha^\beta=\mua-\beta$, and \eqref{torpar} turns into
$$
\begin{aligned}
    \mathscr{T}(\csys(C),Z_\alpha^\beta)& =\tab \\
    \mathscr{F}(\csys(C),Z_\alpha^\beta)& =\fab
    \end{aligned}
$$
and hence \eqref{eq:tilt} reads
$$
\csys(C)^{\dagger}=\csysab(C)
$$

Now we will build a Bogomolov-Gieseker inequality for $(\csys(C),Z_\alpha^\beta)$.

\begin{Lemma}\label{quadratic}
Let $\mathcal{E}$ be a coherent system of type $(n,d,k)$. If $\cale$ is $\mu_\alpha^\beta$-semistable, then 
\begin{equation} \label{equDLT}
\Delta(\mathcal{E})=k(d+n-k)+pd^2+qn^2+uk^2 \geq 0
\end{equation}
for any $p, q \geq 0$ and $u\geq 1$.
\end{Lemma}

\begin{proof}
As $\mu_\alpha^\beta=\mua-\beta$, it follows that $\mathcal{E}$ is $\mu_\alpha$-semistable. If $k=0$, we are done. So assume $k>0$. If $n > 0$, then $d + n-k\geq 0$ by Lemma \ref{cota}, and hence $\Delta(\mathcal{E})\geq 0$.  Now if $n=0$, then $d\geq 0$, and hence
     \[\begin{aligned}
     \Delta(\mathcal{E}) & =k(d-k)+pd^2+uk^2 \\
     &= kd+pd^2+(u-1)k^2 \geq 0 \\
     \end{aligned}\]
and we are done again.  
\end{proof}

The following result provides the expression for $Z^\dagger$ in our context, and in particular proves that the pre-stability conditions presented in Section \ref{sec:stab-tilt} satisfy the support property at least for certain values of parameters.

\begin{Theorem}\label{supproperty}
Let $\alpha,\beta\in\mathbb{Q}_{>0}$, $q,u,t,A\in\mathbb{R}$, with $q\geq 0$, $u\geq 1$, $t>0$ and $\alpha\geq u-1$. Consider the additive homomorphism $\big(Z_\alpha^\beta\big)^\dagger_t:\Lambda =\mathbb{Z}^3\to\mathbb{C}$ defined by
\begin{equation} \label{equztt}   
\big(Z_\alpha^\beta\big)_t^\dagger = \big(1+t(q-A\beta)\big)n+t\big(A-p\beta\big)d+t\big(1+A\alpha-(p\alpha+1)\beta\big)k+\sqrt{-1}(d+\alpha k-\beta n),
\end{equation}
where $p:=(\alpha-u+1)/\alpha$. Then $\tau_\alpha^{\beta,\gamma}=(\csys_\alpha^\beta(C),(Z_\alpha^\beta)_t^\dagger)$ is a stability condition in $D^b(\csys(C))$.

In particular, if $\gamma\in\mathbb{Q}$, with $\gamma\geq(\beta^2+2\alpha\beta-\alpha)/(\alpha(\beta-1))$ if $\beta>1$ or $\gamma > 1+(\beta(1-\alpha))/2\alpha$ if $\beta<1$, then $(\csys_\alpha^\beta(C),Z_\alpha^{\beta,\gamma})$ is a stability condition in $D^b(\csys(C))$.
\end{Theorem}

\begin{proof}
For $(n,d,k)\in\mathbb{Z}^3$, let $\Delta:=k(d+n-k)+pd^2+qn^2+uk^2$ be as in \eqref{equDLT}. Let also $\Re(Z_\alpha^{\beta}):=-(d+\alpha k-\beta n)$ and $\Im(Z_\alpha^{\beta}):=n$ be as in \eqref{equnab}. We may write
\begin{equation}
\label{equDBG}
\Delta =\Re(Z_\alpha^{\beta})\,\Delta_R + \Im(Z_\alpha^{\beta})\,\Delta_I
\end{equation}
as in \eqref{bg-ineq}, where $\Delta_R=An+Bd+Ck$ for  $A,B,C \in \mathbb{R}$ and  $\Delta_I=Dn+Ed+Fk$ for $D,E,F \in \mathbb{R}$. Then \eqref{equDBG} yields the system of equations 
$$ \begin{array}{lll}
  A\beta +D=q & -B = p & -C\alpha =u-1 \\
    -A+B\beta+E =0\ \ \ \ & -A\alpha +C\beta +F =1\ \ \ \  & -C-B\alpha =1   
\end{array} $$
which yields the equalities
$$
B=-p\ \ \ \ \   C=p\alpha-1\ \ \ \ \ D=q-A\beta\ \ \ \ \ E = A+p\beta\ \ \ \ \ F = 1+A\alpha-(p\alpha-1)\beta 
$$
and 
\begin{equation}
\label{App0}
p = \frac{1}{\alpha}\left(\frac{\alpha-u+1}{\alpha}\right)
\end{equation}
and $p\geq 0$ since $\alpha\geq u-1$. So, following \eqref{eq:zdagger}, we set
\[\begin{aligned}
    \big(Z_\alpha^\beta\big)_t^\dagger &:=  \Im(Z_\alpha^\beta)+t\Delta_I-\sqrt{-1}\,\Re(Z_\alpha^\beta)\\
    &= n(1+t(q-A\beta))+t(A-p\beta)d+t(1+A\alpha-(p\alpha+1)\beta)k+\sqrt{-1}(d+\alpha k-\beta n)
\end{aligned}\] 
and the statements, based on \cite[Section 2]{PT}, developed along this section yield the first claim, that is, $ (\csys_\alpha^\beta(C),\big(Z_\alpha^\beta\big)_t^\dagger)$ is a stability condition.

To prove the second, comparing $\Re((Z_\alpha^\beta\big)_t^\dagger)$ and $\Re(Z_\alpha^{\beta, \gamma})$, we get
\[n(1+t(q-A\beta))+t(A-p\beta)d+t(1+A\alpha-(p\alpha+1)\beta)k =  \gamma n+d-k.\]
This gives the following equations
\begin{align}
\label{App1}    t(q-A\beta)&=\gamma-1 \\
\label{App2}    t(A+p\beta)&=1 \\
\label{App3}    t\big((p\alpha-1)\beta-1-A\alpha\big)&=1
\end{align}
As $p=p(\alpha,u)$ as in \eqref{App0}, we have to find $q,u,t,A$, satisfying the equations above and also with $q\geq 0$, $0\leq u-1\leq\alpha$ and $t>0$. Now \eqref{App2} yields 
\begin{equation}
\label{App4}
A=\frac{1}{t}-p\beta.
\end{equation} 
Replacing \eqref{App4} in \eqref{App3} we get
\begin{equation}
\label{App5}
    2p\alpha\beta = \frac{\alpha+1}{t}+\beta+1.
\end{equation}
Replacing \eqref{App0} in \eqref{App5} yields
\begin{equation}
\label{App6}
u=\frac{t\beta(\alpha+2)-\alpha(\alpha+t+1)}{2t\beta}.
\end{equation}
Now, combining \eqref{App4} and \eqref{App5} we get
\begin{equation}
\label{App7}
A = \frac{\alpha-t(\beta+1)-1}{2t\alpha} 
\end{equation}
while \eqref{App1} and \eqref{App7} produces
\begin{equation}
\label{App8}
q = \frac{2\alpha(\gamma-1)+\beta(\alpha-t(\beta+1)-1)}{2t\alpha}.
\end{equation}
So we  found $q,u,t$ and $A$ satisfying \eqref{App1}-\eqref{App3}, with $t>0$ being arbitrary. Now let us check when $q\geq 0$ and $0\leq u-1\leq\alpha$. The first inequality yields
\begin{equation}
\label{App9}
t\leq\frac{2\alpha(\gamma-1)+\beta(\alpha-1)}{\beta(\beta+1)}.
\end{equation}
Now note that $u\leq\alpha+1$ is always satisfied. On the other hand, $u\geq 1$ reads
\begin{equation}
\label{App10}
t(\beta-1)\geq\alpha+1
\end{equation}
If $\beta>1$, then, combining \eqref{App9} and \eqref{App10}, we can find $t$ if and only if
$$
\gamma\geq\frac{\beta^2+2\alpha\beta-\alpha}{\alpha(\beta-1)}
$$
If $\beta<1$, then, by \eqref{App10}, $t\leq (\alpha+1)/(1-\beta)$, which is positive. So we can find $t$ if and only if the numerator on the right-hand side of \eqref{App9} is positive as well. This yields
$$
\gamma > 1+\frac{\beta(1-\alpha)}{2\alpha}
$$
Finally, if $\beta=1$, then \eqref{App10} precludes the existence of such a $t$.
\end{proof}

Let us now consider the set
$\mathsf{PS} := \R_{\ge0}\times\R_{\ge0}\times\R_{\gamma>1}$, and 
the subset $\mathsf{S}\subset\mathsf{PS}$ such that $\beta\ne 0,1$ and either
\begin{itemize}
\item [(i)]  $\gamma\geq(\beta^2+2\alpha\beta-\alpha)/(\alpha(\beta-1))$, if $\beta >1$,
\item [(ii)] or $\gamma > 1+(\beta(1-\alpha))/2\alpha$, if $\beta<1$.
\end{itemize}

We already showed in the end of Section \ref{sec:stab-tilt} that if $(\alpha,\beta,\gamma)\in\mathsf{PS}\cap\Q^3$, then \linebreak $\tau_\alpha^{\beta,\gamma}\in\stablf(D^b(\csys(C)))$.   Moreover, we conclude from Theorem \ref{supproperty} that if $(\alpha, \beta, \gamma) \in\mathsf{PS}\cap\Q^3$ then $\tau_\alpha^{\beta,\gamma}\in\stab(D^b(\csys(C)))$.

Our next statement follows directly from Theorem \ref{supproperty} and \cite[Corollary 2.22]{PT}.

\begin{Corollary} \label{cor:map2stab}
The map $\mathsf{S} \to \stab(D^b(\csys(C)))$ defined by $(\alpha,\beta,\gamma)\mapsto \tau_\alpha^{\beta,\gamma}$ is continuous.
\end{Corollary}

To conclude this section, we will compare the tilted (pre-)stability conditions \linebreak $\tau_\alpha^{\beta,\gamma}=(\csys_\alpha^\beta(C),Z_\alpha^{\beta,\gamma})$ with the standard ones $\sigma_\alpha=(\csys,Z_\alpha)$ constructed in Section \ref{sec:stab-cs} and among themselves.

\begin{Lemma} \label{action2}
For any $\alpha,\alpha'>0$, $\beta\ge0$ and $\gamma>1$, $\tau_\alpha^{\beta,\gamma}\not\in\widetilde{GL}^+(2,\R)\cdot\sigma_\alpha$.
\end{Lemma}

This means that the (pre-)stability conditions constructed via tilting are never equivalent to a standard stability condition.

\begin{proof}
Let $T= \begin{pmatrix} A & B \\ C & D \end{pmatrix} \in GL^+(2,\mathbb{R})$.
We argue that there is no $T\in GL^+(2,\R)$ such that $T\circ Z_{\alpha'} = Z_\alpha^{\beta,\gamma}$, which implies that $Z_\alpha^{\beta,\gamma}$ cannot be in the $\widetilde{GL}^+(2,\R)$ orbit of $\sigma_\alpha$.

Indeed, the equality $T\circ Z_{\alpha'} = Z_\alpha^{\beta,\gamma}$ becomes
$$ \begin{pmatrix} A & B \\ C & D \end{pmatrix}
\begin{pmatrix} -d-\alpha' k \\ n \end{pmatrix} =
\begin{pmatrix} d+\gamma n-k \\ d+\alpha k- \beta n 
\end{pmatrix} $$
which is equivalent to 
\[\begin{aligned}
    -A(d+\alpha' k)+Bn &=d+\gamma n-k, ~~{\rm and} \\
    -C(d+\alpha' k)+Dn &=d+\alpha k- \beta n.
\end{aligned}\]
It follows that $B=\gamma$, $C=-1$, $D=-\beta$, and  $A=-1$. In addition, we would have that  $\alpha'k=-k$, which is impossible because $\alpha'\ge0$.
\end{proof}

\begin{Lemma} \label{action3}
Let $\alpha,\beta,\gamma$ and $\alpha',\beta',\gamma'$ be such that $\tau_\alpha^{\beta,\gamma}$ and $\tau_{\alpha'}^{\beta',\gamma'}$ are pre-stability conditions on $D^b(\csys(C))$. Then:
\begin{itemize}
\item [(i)] If $\gamma \neq \gamma'$, then $\tau_\alpha^{\beta,\gamma}\not\in\widetilde{GL}^+(2,\R)\cdot\tau_{\alpha'}^{\beta',\gamma'}$;
\item [(ii)] If $(\alpha-\alpha')\gamma-(\alpha'+1)\beta \neq -\beta'(\alpha+1)$, then
$\tau_\alpha^{\beta,\gamma}\not\in\widetilde{GL}^+(2,\R)\cdot\tau_{\alpha'}^{\beta',\gamma'}$.
\end{itemize}
\end{Lemma}

\begin{proof}
Let $T= \begin{pmatrix} A & B \\ C & D \end{pmatrix} \in GL^+(2,\mathbb{R})$. If $T \circ Z_\alpha^{\beta,\gamma}=Z_{\alpha'}^{\beta', \gamma'}$, then
$$ \begin{array}{cc}
\vspace{2mm} A+B = 1 \hspace{5mm} & \hspace{5mm} C+D =1\\ 
\vspace{2mm}    -A+\alpha B = -1\hspace{5mm} & \hspace{5mm} C\gamma-D\beta = -\beta' \\
A\gamma-B\beta=\gamma' \hspace{5mm} & \hspace{5mm} \alpha D-C =\alpha' \\    
\end{array} $$
which yields $A=1$, $B=0$, $C=\frac{\alpha-\alpha'}{\alpha+1}$, $D= \frac{\alpha'+1}{\alpha+1}$; in addition, we must also have that
$$ \gamma=\gamma' ~~{\rm and}~~ (\alpha-\alpha')\gamma-(\alpha'+1)\beta=-\beta'(\alpha+1). $$
Therefore, if neither of these equalities holds, then the $\widetilde{GL}^+(2,\mathbb{R})$ orbit of $\tau_\alpha^{\beta,\gamma}$ does not contain $\tau_{\alpha'}^{\beta',\gamma'}$.
\end{proof}


\section{Semistability of Tilted Torsion Systems} \label{sec:sst-torsion}

In this section, we analyse the $\mu_\alpha^{\beta,\gamma}$-semistability of certain objects in $\csys^\beta_\alpha(C)$ and explore related properties. Our focus is on objects that arise from torsion sheaves.

Following \cite[Exercise 6.9]{Macri-Schmidt}, say that an object in an abelian category $\mathscr{A}$ is \emph{minimal} if it does not have any nontrivial subobject (or quotient)\footnote{the term \emph{simple} is also frequently used in the context of abelian categories}.  Denote by $\Min(\mathscr{A})$ the set of minimal objects in $\mathscr{A}$. Given $E\in\mathscr{A}$, say $E'\in\mathscr{A}$ is an \emph{elementary transformation of $E$ along $S$}, if there is an exact sequence 
    \begin{equation}\label{ElmT}
        0\to E' \to E \to S \to 0.
    \end{equation}  
where $S\in \Min(\mathcal{A})$.

\begin{Lemma}
\label{HNFCAT}
Let $E$ be an object in $\mathscr{A}$ which is not semistable and with Harder--Narasimhan filtration $0 =E_0 \subset E_1 \subset \ldots \subset E_l=E$. Let $E'$ be an elementary transformation of $E$ along $S$. Consider the morphisms $\psi_i:E_i\to E\to S$ and set $E_i':=\ker(\psi_i)$. If $\psi_1\ne0$, then
$0\subset E_1' \subset\cdots\subset E_l'=E'$ is the Harder--Narasimhan filtration of $E'$, and also $E'_i/E_{i-1}'=E_i/E_{i-1}$ for all $i\in\{1,\ldots,l\}$. In particular, $E'$ is not semistable.
 \end{Lemma}
 
 \begin{proof}
First note that if $\psi_i\neq0$, then it is surjective. Indeed, as $\im(\psi_i)\subset S$ and $S$ is minimal, it follows that $\im(\psi_i)=S$, that is, $\psi_i$ is surjective. Now if $\psi_1\neq 0$, then, by construction, all $\psi_i$ do not vanish, and hence, all are surjective. 

Therefore, for each $i=1,\dots,k$ we have the diagram
$$ \xymatrix{
& 0 \ar[d] & 0 \ar[d] & & \\
0 \ar[r] & E_{i-1}' \ar[r]\ar[d] & E_{i-1} \ar[r]^{\psi_{i-1}}\ar[d] & S \ar[r]\ar@{=}[d] & 0 \\
0 \ar[r] & E_{i}' \ar[r]\ar[d] & E_i \ar[r]^{\psi_{i}}\ar[d] & S \ar[r] & 0 \\
& E_{i}'/E_{i-1}' \ar@{=}[r]\ar[d] & E_{i}/E_{i-1}\ar[d] & & \\
& 0  & 0 & &
} $$
It follows that $0\subset E_1' \subset\cdots\subset E_l'=E'$ is a filtration whose factors are semistable and $M(E_{i}'/E_{i-1}')>M(E_{i+1}'/E_{i}')$. Therefore, it is the Harder--Narasimhan filtration of $E'$. Note that we needed $E$ not to be semistable to guarantee the existence of the diagram above, with two different (horizontal) sequences with the same quotient $S$. In particular, the Harder--Narasimhan filtration of $E'$ is of the same length as that of $E$, which is at least $2$, so $E'$ is not semistable either.
\end{proof}

Now given a point $P\in C$, we denote the structure sheaf of $P$, viewed as a closed subscheme of $C$, by $\mathcal{O}_P$. Note that 
\begin{equation}
\label{MIN}
\Min(\csys(C))=\{\{(\mathcal{O}_P,0)\}_{P\in C}, \{(0,V_1)\}_{\dim(V)=1}\}.
\end{equation}
For simplicity, we will denote by $V_1$ any $k$-vector space of dimension $1$.

\begin{Definition}
\emph{Let $\eee\in \csys(C)$ and let $\mathcal{E}'$ be an elementary transformation of $\mathcal{E}$ fitting into an exact sequence 
    \begin{equation}\label{ElmType}
        0\to \mathcal{E}' \to \mathcal{E} \to \mathcal{S} \to 0
    \end{equation}  
where $\mathcal{S}\in \Min(\csys(C))$. We say that $\mathcal{E}'$ is (an elementary transformation of $\eee$) of \emph{type I} if $\mathcal{S}=(\mathcal{O}_P,0)$ for some $P\in C$, and it is of \emph{type II} if $\mathcal{S}=(0,V_1)$ for some $V_1$.}
\end{Definition}

\begin{Remark}
\label{ElemTypeII}
\emph{Let $\eee=(E,V,\varphi)\in\csys(C)$ and $\eee'=(E',V',\varphi')$ be an elementary transformation of $\eee$ of type I. If $E$ and $E'$ are locally free, then $E'$ is an elementary transformation of $E$ in the sense of (cf. \cite[p.~5]{BramLange}) or a \emph{Hecke modification} (cf. \cite{NarRam}). Here, we distinguish between types I and II of elementary transformations since our categorical approach to coherent systems involves two types of minimal objects. Note that a surjective morphism $E \to \mathcal{O}_P$ may not induce an elementary transformation of type I. In fact, the existence of any such is closely related to the notion of base points of a coherent system, which we introduce right away.}
\emph{On the other hand, if $\mathcal{E}$ is of type $(n,d,k)$, with $k>0$, one can always find an elementary transformation of type II taking $\mathcal{E}':=(E,V', \varphi|_{V'})$ where $V'\subset V$ and $\dim V'= \dim V-1$.}
\end{Remark}

\begin{Definition}
\label{defpb}
\emph{Let $\eee \in \csys(C)$. A point $P \in C$ is a \emph{base point} of $\eee$ if for all $x\in V$ 
we have that $x$ vanishes at $P$, i.e., $x\in \mmp E_P$.  A coherent system is said to be \emph{base point free} if it has no base points.}
\end{Definition}

\begin{Remark}
\emph{The definition above works for our purposes here, and matches the usual definition when $C$ is smooth. However, it doesn't recover the case where $C$ is singular and the system is \emph{linear}, i.e., $\rank(E)=1$. In fact, given a linear system $\mathcal{L}=(L,V)$, following \cite[p.~8]{KM2} (based on \cite[p.~198]{RSt}), a point $P\in C$ is said to be a base point of $\mathcal{L}$ if  for all $x\in V$, the injection $\varphi_{P,x}:\oo_P\to L_P$ is not surjective. Note that if $L$ is invertible (which always holds if $C$ is smooth), then $\mathcal{L}$ is base point free if and only if $L$ is generated by $V$. If $L$ is torsion-free of rank $1$ but not invertible, then $\mathcal{L}$ can be generated by $V$ and, even that, admit a base point, which is then called \emph{irremovable} (\cite[p.~198]{RSt}). The motivation for this definition has some geometric content. It formally describes, for instance, a curve $C$ lying on a cone $S$ passing through the vertex $P$, such that $P$ is a singularity of $C$. Then the linear system $\mathcal{L}$ cut out by the rulings of $S$ is such that $P$ is an irremovable base point of $\mathcal{L}$.} 
\end{Remark}

\begin{Proposition}
\label{TEPB}
    Let $\eee=(E,V) \in \csys(C)$ and $P\in{\rm Supp}(E)$. Then there exists a surjective morphism $\eee \to (\mathcal{O}_P,0)$ if and only if $P$ is a base point of $\eee$. 
\end{Proposition}

\begin{proof}
As $P\in{\rm Supp}(E)$, the fibre $E_{(P)}:=E_P/\mmp E_P$ does not vanish. So take a linear surjective map $h_P:E_{(P)}\to k_P=H^0(\op)$. It lifts to a surjective map of sheaves $h:E\to \op$ and conversely.
And we have the diagram
\begin{equation}
\xymatrix{
    V \ar[r] \ar[d] & H^0(E) \ar[d]^{H^0(h)}\\
    0 \ar[r] & H^0(\mathcal{O}_P)}
\end{equation} 
which is commutative if and only if $P$ is a base point of $\eee$.
\end{proof}

Let $P\in C$ be any point. Consider the torsion-free sheaf $M_{\{P\}}$ of rank $1$ on $C$ defined by the exact sequence
$$
0\to M_{\{P\}}\to \oc \to \oo_P\to 0
$$
Given a coherent sheaf $E$, for simplicity, we will write $E(-P):=M_{\{P\}}E$. Note, for instance, that if $\mathcal{E}=(E,V)$ is a system with $E$ invertible, and with a base point $P$, then $(E(-P),V)$ is an elementary transformation of $\mathcal{E}$ of type I. Note that if $P$ is singular, then $E(-P)$ is not invertible, even when defining an elementary transformation arising from an invertible sheaf. Similar examples can be obtained for higher rank systems.

\begin{Lemma}\label{Notunsta}
Let $\mathcal{E}\in\csys(C)$. If $\eee \in \mathscr{T}_\alpha^\beta$ (resp. $\fab$) and admits an elementary transformation $\mathcal{E}'$ which is in $\mathscr{F}_\alpha^\beta$ (resp. $\tab$), then $\eee$ is $\mu_\alpha$-semistable.
\end{Lemma}

\begin{proof}
Let $0 =\eee_0 \subset \eee_1 \subset \ldots \subset \eee_l=\mathcal{E}$ be the Harder--Narasimhan filtration of $\eee$ and $0\subset \eee_1' \subset\cdots\subset \eee_{p}'=\eee'$ be the one of $\eee'$. Assume $\eee$ is not $\mu_\alpha$ semistable. Then Lemma \ref{HNFCAT} implies $l=p$ and $\eee _1=\eee_1'$. As $\eee_1'\subset \eee'$ and $\eee' \in \fab$, then $\mu_\alpha(\eee_1')\leq\beta$. On the other hand, $\eee_1$ destabilizes $\eee$, which is in $\tab$, so $\mu_\alpha(\eee_1')>\mua(\eee')>\beta$. As $\eee_1=\eee_1'$, we get a contradiction. A similar argument holds for the respective part.
\end{proof}

For the next result, note that, by \eqref{MIN}, any minimal object $\eee\in\csys(C)$ has rank $0$. Thus if $\eee\twoheadrightarrow\mathcal{Q}$, then  $\mathcal{Q}$ has rank $0$ too. Therefore $\mu_\alpha(\mathcal{Q})=\infty >\beta$ and so $\eee\in\mathcal{T}_\alpha^\beta$. Thus $\eee$ is naturally seen as an object in $\csysab(C)$, with $\hhz(\eee)=\eee$ and $\hhm(\eee)=0$.

\begin{Lemma} 
\label{prpozz1}
Let $\eee\in \Min(\csys(C))$. Then the following holds.
\begin{itemize}
\item[(i)] $\fff$ is a nonzero proper subobject of $\eee$ in $\csysab(C)$ if and only if $\fff\in\tab$ and fits into an exact sequence in $\csys(C)$, 
\begin{equation}
\label{equrs1}
          0  \to \mathcal{F}' \to 
       \mathcal{F} \to \eee \to  0
\end{equation} 
with $\mathcal{F}'\in\mathcal{F}_\alpha^\beta$. Also, if $\mathcal{F}$ is of type $(n,d,k)$, then $n\neq 0$ and
\begin{itemize}
    \item [(a)]If $\mathcal{E}=(\oo_P,0)$, then 
\begin{equation}
\label{equrs2}
d-1+\alpha k \leq n\beta < d+\alpha k; 
\end{equation}
\item [(b)] If $\mathcal{E}=(0,V_1)$,  then 
\begin{equation}
\label{equrs21}
d +\alpha (k-1) \leq n\beta < d+\alpha k.
\end{equation}
\end{itemize}
\item[(ii)] $\eee=(\oo_P,0)$ is $\mu_\alpha^{\beta,\gamma}$-semistable (resp. stable) if and only iffor every $\fff$ as in (i) we have
\begin{equation}
\label{equrs3}
 k(\alpha+1) \leq n(\beta+\gamma)   \ \ (resp.\,<)
\end{equation}
\item[(iii)] $\eee=(0,V_1)$ is $\mu_\alpha^{\beta,\gamma}$-semistable (resp. stable) if and only if for every $\fff$ as in (i) we have
\begin{equation}
\label{equrs31}
n(\beta-\alpha\gamma) \leq d(\alpha+1)\ \ (resp.\,<)
\end{equation}
\end{itemize}
\end{Lemma}

\begin{proof} 
To prove (i), let $0\neq\mathcal{F}\subsetneq\eee$. Consider the exact sequence 
   \begin{equation}
     0\to \mathcal{F} \to \eee \to \mathcal{G} \to 0
    \end{equation}
in $\csys_\alpha^\beta(C)$. It yields the following exact sequence 
  \begin{equation}
  \label{equlex}
       0 \to \hhm(\mathcal{F}) \to  0 \to  \hhm(\mathcal{G}) \to 
         \hhz(\mathcal{F}) \to \eee \to \hhz(\mathcal{G}) \to 0
   \end{equation}
    in $\csys(C)$.
Thus 
$\hhm(\mathcal{F})=0$, and, as $\eee$ is minimal, 
either $\eee \cong \hhz(\mathcal{G})$, or $\hhz(\mathcal{G})=0$.
    If $\eee\cong \hhz(\mathcal{G})$, then $\hhm(\mathcal{G}) \cong \hhz(\mathcal{F})$. Now  $\tab\cap\fab=0$, so $\hhz(\mathcal{F})=0$. But $\hhm(\mathcal{F})=0$ too, thus $\mathcal{F}=0$, which is precluded as $\mathcal{F}$ is assumed to be nonzero. Therefore $\hhz(\mathcal{G})=0$. So we conclude that $\mathcal{F}\cong \hhz(\mathcal{F})\in\tab$, and $\mathcal{G}\cong \hhm(\mathcal{G})[1]$ with $\hhm(\mathcal{G})\in\fab$. Set  $\mathcal{F}':=\hhm(\mathcal{G})$. Then \eqref{equlex} turns into 
\begin{equation}
    \label{equrst}
          0  \to \mathcal{F}' \to 
       \mathcal{F} \to \eee \to  0
   \end{equation}
which agrees with \eqref{equrs1}. Conversely, \eqref{equrst} yields
$$  0\to  \mathcal{F} \to \eee \to \mathcal{F}'[1]  \to 0
$$
which is an exact sequence in $\csys_\alpha^\beta(C)$, and hence $\mathcal{F}$ is a nonzero subobject of $\eee$. This proves the first statement of (i). To prove the second, say $\mathcal{F}$ is of type $(n,d,k)$ and $\mathcal{F}'$ of type $(n',d',k')$. Then, for (a), $n'=n$, $d'=d-1$, $k'=k$ and, and for (b), $n'=n$, $d'=d$, and $k'=k-1$. As $\mathcal{F}'\in\fab$, then $n'>0$. Also, $\mathcal{F}\in\tab$. Thus, for (a), we get
$$
\mu_\alpha(\mathcal{F}')=\frac{(d-1)+\alpha k}{n}\leq\beta<\frac{d+\alpha k}{n}=\mu_\alpha(\mathcal{F})
$$
while for (b) we get
$$
\mu_\alpha(\mathcal{F}')=\frac{d+\alpha (k-1)}{n}\leq\beta<\frac{d+\alpha k}{n}=\mu_\alpha(\mathcal{F})
$$
which yields \eqref{equrs2} and \eqref{equrs21}.

To prove (ii), note that (i) characterizes the subobjects $\fff\subset\eee=(\oo_P,0)$. So the latter is $\mu_\alpha^\beta$-semistable if and only if for any such $\fff$,  we have 
$$
\mu_\alpha^{\beta,\gamma}(\mathcal{F})=-\frac{d+\gamma n-k}{d+\alpha k-\beta n}\leq  -1=\mu_\alpha^{\beta,\gamma}((\oo_P,0))
$$
which yields \eqref{equrs3} since $d+\alpha k-\beta n>0$ because $\mathcal{F} \in \tab$. The proof of (iii) is similar.
\end{proof}

A stability condition on the derived category of coherent sheaves on a scheme is called \textit{geometric} if all the skyscraper sheaves $\oo_P$ (which are the minimal objects in a category of sheaves) are stable of the same phase (see \cite[Definition 2.1]{BM11}). Following this line, we say that a stability condition on $D^b(\csys(C))$ is \textit{geometric} if the objects $(\oo_P,0)$ and $(0,V_1)$ are stable. 
The following result states the stability of $(\oo_P,0)$ and $(0,V_1)$ in $\csys_\alpha^\beta(C)$.

\begin{Theorem}\label{geometric}
For every $(\alpha,\beta,\gamma)\in\mathsf{PS}$, the following holds:
\begin{itemize}
\item[(i)]
$(\oo_P,0)$ is $\mu_\alpha^{\beta,\gamma}$-stable; 
\item[(ii)] $(0,V_1)$,  is $\mu_\alpha^{\beta,\gamma}$-stable.
\end{itemize}
\end{Theorem}
\begin{proof}
To prove (i), by Lemma \ref{prpozz1}.(i), any $0\neq\fff\subsetneq \eee=(\oo_P,0)$ fits into an exact sequence like \eqref{equrs1}. Also, $\fff\in\tab$ and if $\fff$ is of type $(n,d,k)$, then $n>0$ and
\begin{equation}
\label{equrs41}
d+\alpha k-1\leq\beta n<d+\alpha k
\end{equation}
Besides, if $\fff$ destabilizes $(\oo_P,0)$, then, by Lemma \ref{prpozz1}.(ii), 
\begin{equation}
\label{equrs51}
 n(\beta+\gamma) \leq (\alpha+1)k  
\end{equation}
Combining \eqref{equrs41} and \eqref{equrs51}, we get
$$
d-1+\alpha k \leq \beta n < (\alpha+1) k - \gamma n<\alpha k +k - n
$$
which yields 
\begin{equation}
\label{descont1}
k\geq d+n
\end{equation}
On the other hand, by \eqref{equrst}, $\fff'\in\fab$ is an elementary transformation of $\fff\in\tab$. Thus Lemma \ref{Notunsta} implies $\mathcal{F}$ is $\mu_\alpha$-semistable. Now $n,\beta >0$ and $\gamma>1$, so (\ref{equrs51}) yields $k>0$, and Theorem \ref{cota}.(i) yields.
\begin{equation}
\label{descont2}
k\leq d+n
\end{equation}
So from \eqref{descont1} and \eqref{descont2} we get $k=d+n$, that is, equality holds in Theorem \ref{cota}.(i). But equality holds in Theorem \ref{cota}.(i) if and only if $\mathcal{F}=|\mathcal{O}_C^{\oplus n}|$ or $\fff= |\oplus_{i=1}^{n}\mathcal{O}_{\mathbb{P}^1}(a_i)|$  and $a_i\geq 0$. As both systems have no base points, then $\fff$ cannot admit an elementary transformation owing to Proposition \ref{TEPB}, a contradiction. So $(\oo_P,0)$ is $\mu_\alpha^{\beta,\gamma}$-stable.

To prove (ii), proceeding as above, by Lemma \ref{prpozz1}. (ii), a destabilizer $\fff\in\tab$ of type $(n,d,k)$ would satisfy 
\begin{equation}
\label{equrs61}
d+\alpha (k-1) \leq n\beta <d+\alpha k, 
\end{equation}
\begin{equation}
\label{equrs71}
d(\alpha+1) \leq n(\beta-\alpha\gamma).
\end{equation}
Now from \eqref{equrs61} and \eqref{equrs71}, it follows that $d(\alpha+1)+n\alpha\gamma < n\beta <d+\alpha k$ which implies $d+n<k$.  But $d\geq 0$ whenever $k=0$ because $n\beta < d+\alpha k$. Also, $\fff$ is $\mua$-semistable by Lemma \ref{Notunsta}, which contradicts Theorem \ref{cota}.(i). So $(0,V_1)$ is $\mu_\alpha^{\beta,\gamma}$-stable. 
\end{proof}

\begin{Corollary} \label{opcomp}
For every $(\alpha,\beta,\gamma)\in\mathsf{PS}$ and every point $P\in C$, the system $|\oo_P|$ is $\mu_\alpha^{\beta,\gamma}$-stable.
\end{Corollary}

\begin{proof}  
Let $0\neq\mathcal{F}\subsetneq\eee$ be a nontrivial, proper subobject, and consider the exact sequence 
\begin{equation}
0\to \mathcal{F} \to |\oo_P| \to \mathcal{G} \to 0
\end{equation}
in $\csys_\alpha^\beta(C)$. It yields the following exact sequence 
\begin{equation}
0 \to \hhm(\mathcal{F}) \to  0 \to  \hhm(\mathcal{G}) \to 
\hhz(\mathcal{F}) \to |\oo_P| \stackrel{f}\to \hhz(\mathcal{G}) \to 0
\end{equation}
in $\csys(C)$. Thus $\hhm(\mathcal{F})=0$ and $\hh^0(\mathcal{F})=\fff$. Set $\mathcal{K}:={\rm ker}f$; since this is a subsystem of $|\oo_P|$ there are three possibilities: either $\mathcal{K}=0$, or $\mathcal{K}=|\oo_P|$, or $\mathcal{K}=(\oo_P,0)$, which is the only proper, nontrivial subobject of $|\oo_P|$.

If $\mathcal{K}=0$, then $\hhm(\mathcal{G})\cong\hh^0(\mathcal{F})$ which yield $\hh^0(\fff)=0$ and $\fff=0$ contradicting our initial assumption. 

If $\mathcal{K}=(\mathcal{O}_P,0)$, then $\hhz(\mathcal{G})=(0,V_1)$ and we obtain the following exact sequences in $\csys(C)$:
$$ 0 \to  \hhm(\mathcal{G}) \to  \fff \to (\mathcal{O}_P,0) \to 0 ~~ {\rm and}$$
$$ 0 \to (\mathcal{O}_P,0) \to |\oo_P| \to (0,V_1) \to 0. $$
The first of these induces the following exact sequence in $\csysab(C)$:
$$ 0 \to \fff \to (\mathcal{O}_P,0) \to \hhm(\mathcal{G})[1] \to 0, $$
so $\calf$ is a subobject of $(\mathcal{O}_P,0)$ in $\csysab(C)$. It follows that 
\[\mu_\alpha^{\beta, \gamma}(\fff) \leq \mu_\alpha^{\beta, \gamma}(\oo_P,0) = -1 < 0=\mu_\alpha^{\beta, \gamma}(|\oo_P|), \]
so $\calf$ does not destabilizes $|\oo_P|$.

Finally, if $\mathcal{K}=|\mathcal{O}_P|$, then $\hhz(\mathcal{G})=0$, so $\calg=\hhm(\mathcal{G})[1]$, and we get the following exact sequence in  $\csys(C)$:
\[ 0 \to \calg \to  \fff \to |\mathcal{O}_P| \to 0 \]
with $\calg\in\fab$ and $\calf\in\tab$. In addition, if $v(\calf)=(n,d,k)$, then $v(\calg)=(n,d-1,k-1)$, $n>0$. Setting $\calf=(F,U,\phi)$ and $\calg=(G,W,\psi)$, the latter sequence can be unpacked into the following commutative diagram
$$ \xymatrix{
0 \ar[r] & W \ar[r] \ar[d]^{\psi} & U \ar[r] \ar[d]^{\phi} & V_1 \ar[r] \ar[d]^{\simeq} & 0 \ar[d] \\ 
0 \ar[r] & H^0(G) \ar[r] & H^0(F) \ar[r] & H^0(\oo_P) \ar[r] & H^1(E)
}$$
Note that $\psi$ is injective because $\calg\in\fab$; it follows that $\phi$ is also injective, and we can use Proposition \ref{TorHNF}(ii)(b) to conclude that $k\le d+n$. Hence, we have $\mu_\alpha^{\beta, \gamma}(\fff)<0=\mu_\alpha^{\beta, \gamma}(|\oo_P|)$ so $\fff$ does not destabilizes $|\oo_P|$ and we can conclude that $|\oo_P|$ is $\mu_\alpha^{\beta, \gamma}$-stable. 
\end{proof}


\section{Semistability of Tilted Complete Systems} \label{sec:sst-complete}

In this section, we analyse the $\mu_\alpha^{\beta,\gamma}$-semistability of complete coherent systems (see Definition \ref{DefCohSys}). Also, we determine an actual wall for $\mathcal{E}[1]$ when $\eee$ is $\mu_\alpha$-stable. Our starting point is to study stability properties in the $\gamma$-direction, and from these results, we can apply wall-crossing arguments with respect to the stability function $Z_\alpha^{\beta,\gamma}$ for $\gamma >1$. 

\begin{Lemma}
   Let $\mathcal{E} \in \csys_\alpha^{\beta,\gamma}(C)$ be of type $(n,d,k)$ and $\mathcal{E}'\subset\eee$ be of type $(n',d',k')$. Then $\mu_\alpha^{\beta,\gamma}(\mathcal{E}')-\mu_\alpha^{\beta,\gamma}(\mathcal{E})$ is a linear function in $\gamma$ which is
    \begin{itemize}
        \item[(i)] increasing if $n/(d+\alpha k-\beta n)> n'/(d'+\alpha k'-\beta n')$;
        \item[(ii)] decreasing if $n/(d+\alpha k-\beta n)< n'/(d'+\alpha k'-\beta n')$;
        \item[(iii)] constant if $n/(d+\alpha k-\beta n)= n'/(d'+\alpha k'-\beta n')$.
    \end{itemize}
    In particular, if $\mu_\alpha^{\beta,\gamma'}(\mathcal{E})=\mu_\alpha^{\beta,\gamma'}(\mathcal{E}')$, then
    \begin{itemize}
        \item[(iv)] $\big(\mu_\alpha^{\beta,\gamma}(\mathcal{E}')-\mu_\alpha^{\beta,\gamma}(\mathcal{E})\big)(\gamma-\gamma')>0$, $\forall$ $\gamma\neq\gamma'$ if $n/(d+\alpha k-\beta n)> n'/(d'+\alpha k'-\beta n')$;
        \item[(v)] $\big(\mu_\alpha^{\beta,\gamma}(\mathcal{E}')-\mu_\alpha^{\beta,\gamma}(\mathcal{E})\big)(\gamma-\gamma')<0$, $\forall$ $\gamma\neq\gamma'$ if $n/(d+\alpha k-\beta n)< n'/(d'+\alpha k'-\beta n')$;
        \item[(vi)] $\mu_\alpha^{\beta,\gamma}(\mathcal{E}')-\mu_\alpha^{\beta,\gamma}(\mathcal{E})=0$, $\forall$ $\gamma$ if $n/(d+\alpha k-\beta n)= n'/(d'+\alpha k'-\beta n')$.
    \end{itemize}
    \end{Lemma}

\begin{proof}
    Note that
    \begin{equation} \label{falpha}
    \begin{aligned}
        f(\gamma)&:=\mu_\alpha^{\beta,\gamma}(\mathcal{E}')-\mu_\alpha^{\beta,\gamma}(\mathcal{E}) \\
        &= \gamma \left(\frac{n}{d+\alpha k-\beta n}-\frac{n'}{d'+\alpha k'-\beta n'}\right)+\frac{d-k}{d+\alpha k-\beta n}-\frac{d'-k'}{d'+\alpha k'-\beta n'}.
    \end{aligned}
    \end{equation}
So $f$ is a linear function and items (i)-(iii) follow easily. Items (iv)-(vi) relate to the signal analysis of $f(\gamma)(\gamma-\gamma')$ when $f(\gamma')=0$, and these also follow immediately. \end{proof}

The following result is akin to that in \cite[Lemma 6.2]{BGPMN}, and its proof is analogous.

\begin{Lemma}\label{changesgamma}
    Let $\mathcal{E} \in \csys_\alpha^{\beta,\gamma}(C)$ be of type $(n,d,k)$. Suppose that $\mathcal{E}$ is $\mu_\alpha^{\beta,\gamma}$-stable for $\gamma >\gamma'$, but strictly $\mu_\alpha^{\beta,\gamma}$-semistable for $\gamma=\gamma'$. Then $\mathcal{E}$ is $\mu_\alpha^{\beta,\gamma}$-unstable for all $\gamma < \gamma'$.  
\end{Lemma}

\begin{Corollary}
    Let $\gamma \in  \mathbb{Q}$  and let $\mathcal{E} \in \csys_\alpha^\beta(C)$   of type $(n,d,k)$ such that $\mathcal{E}$ is strictly $\mu_\alpha^{\beta,\gamma}$-semistable.  For any subobject $\mathcal{E}'\subset \mathcal{E}$ of type $(n',d',k')$ such that $\mu_\alpha^{\beta,\gamma}(\mathcal{E}')=\mu_\alpha^{\beta,\gamma}(\mathcal{E})$, the following statements holds:
    \begin{enumerate}
        \item [(i)] If $\mathcal{E}$ is $\mu_\alpha^{\beta,\gamma+\epsilon}$-semistable, then $n/(d+\alpha k-\beta n)> n'/(d'+\alpha k'-\beta n')$.
        \item [(ii)] If $\mathcal{E}$ is $\mu_\alpha^{\beta,\gamma-\epsilon}$-semistable, then $n/(d+\alpha k-\beta n) < n'/(d'+\alpha k'-\beta n')$.
    \end{enumerate}
\end{Corollary}

Following \cite[Section 4]{JM} we define \textit{numerical}, \textit{pseudo} and \textit{actual} walls in the $\gamma$-direction.

\begin{Definition}\label{critical}
Fix $(\alpha,\beta)\in\mathbb{Q}_{\ge0}\times\mathbb{Q}_{\ge0}$.
We say that $\gamma_0\in\R_{>1}$ is
\begin{itemize}
\item[(i)] a \emph{numerical} wall for $(n,d,k) \in \mathbb{Z}^3$, if there is $(n',d',k')\in \mathbb{Z}^3$ non-parallel to $(n,d,k)$ such that
\[-\frac{d'+\gamma_0 n'-k'}{d'+\alpha k'-\beta n'}= -\frac{d+\gamma_0 n-k}{d+\alpha k-\beta n}; \]
\item[(ii)] a \emph{pseudo} wall for $(n,d,k)$ if $\gamma_0$ is a numerical wall and there exists objects $\fff, \eee \in \csys_\alpha^\beta(C)$ with $v(\mathcal{F})=(n',d',k')$, $v(\mathcal{E})=(n,d,k)$ and a nonzero morphism $\calf\to\cale$.
\item[(iii)] an \emph{actual} wall for $(n,d,k)$ if there exists $\eee\in\csys_\alpha^\beta(C)$ of type $(n,d,k)$ 
and a path $\theta:(-1,1) \to \mathbb{R}_{>1}$ such that $\theta(0)=\gamma_0$ and $\eee$ is $\mu_\alpha^{\beta,\theta(t)}$-stable for $t>0$ and $\mu_\alpha^{\beta,\theta(t)}$-unstable for $t <0$.
\end{itemize}
\end{Definition}

Clearly, an actual wall is also pseudo: if $\calf$ is a $\mu_\alpha^{\beta,\theta(t)}$ destabilizing subobject for $t<0$, then $\mu_\alpha^{\beta,\theta(t)}(\calf)<\mu_\alpha^{\beta,\theta(t)}(\cale)$ for $t<0$ with $\mu_\alpha^{\beta,\theta(t)}(\calf)>\mu_\alpha^{\beta,\theta(t)}(\cale)$ for $t>0$ (since $\cale$ is stable in this range), thus $\mu_\alpha^{\beta,\gamma_0}(\calf)=\mu_\alpha^{\beta,\gamma_0}(\cale)$.

The following result describes $\mu_\alpha^{\beta,\gamma}$-semistable objects when $\gamma \gg 0$, akin to the so-called \textit{large volume limit} in the context of stability conditions on the derived category of sheaves on surfaces.

\begin{Lemma} \label{prplim}
Let $\eee\in\csys(C)$ be such that $\rank(\eee)>0$ and $\mu_\alpha(\eee)\neq \beta$. Then
\[\lim_{\gamma \to \infty} \frac{\mu^{\beta, \gamma}_\alpha(\mathcal{E})}{\gamma} = -\frac{1}{\mu_\alpha(\mathcal{E})-\beta}.\] 
\end{Lemma}
\begin{proof}
Say $\eee$ is of type $(n,d,k)$. Then write
$$
        \frac{\mu^{\beta, \gamma}_\alpha(\mathcal{E})}{\gamma}=-\frac{d/\gamma\,+\,n\,-\,k/\gamma}{n(\mu_\alpha(\mathcal{E})-\beta)}
$$
and take the limit to get the desired quotient.
\end{proof}

For the next result, we set 
\[
    \textsf{B}=\textsf{B}_\alpha^\beta := \bigg{\{} \mathcal{E} \in \csys_\alpha^\beta(C) \, \bigg{|}\, \begin{matrix}\ \text{either (i)}\ \hh^{0}(\mathcal{E})=0\ \text{and}\ \hh^{-1}(\mathcal{E})\ \text{is}\ \mu_{\alpha}\text{-semistable}\text{;} \\ \ \ \ \ \, \text{or (ii)}\ \hh^{-1}(\mathcal{E})=0\ \text{and}\  \hh^0(\mathcal{E})\ \text{is}\ \mu_\alpha\text{-semistable}\end{matrix} \bigg{\}}.
\]
With this in mind, we have the following statement.

\begin{Lemma}\label{biggamma}
    If $\mathcal{E} \in \csys_\alpha^\beta(C)$ is $\mu^{\beta,\gamma}_\alpha$-semistable for $\gamma\gg 0$, then $\mathcal{E} \in \textsf{B}$.
\end{Lemma}

\begin{proof}
    We know $\mathcal{E}$ fits into an exact triangle
\begin{equation}
\label{equtri}
0\to \hh^{-1}(\mathcal{E})[1] \to \mathcal{E} \to \hh^{0}(\mathcal{E}) \to 0  
\end{equation}
    where $\hh^{-1}(\mathcal{E}) \in \mathcal{F}^\beta_\alpha$ and $\hh^{0}(\mathcal{E}) \in \mathcal{T}^\beta_\alpha$. 
Now $\mathcal{E}$ is $\mu^{\beta, \gamma}_\alpha$-semistable for $\gamma \gg 0$. Therefore \eqref{equtri} yields
\begin{equation}
\label{equ1<0}
    \frac{\mu^{\beta, \gamma}_\alpha(\hh^{-1}(\mathcal{E}))}{\gamma}=\frac{\mu^{\beta, \gamma}_\alpha(\hh^{-1}(\mathcal{E})[1])}{\gamma} \leq \frac{\mu^{\beta, \gamma}_\alpha(\hhz(\eee))}{\gamma}    
\end{equation}
Taking limits, we get from Proposition \ref{prplim} that
$$
    0 \leq -\frac{1}{\mu_\alpha(\hh^{-1}(\mathcal{E}))-\beta} \leq   -\frac{1}{\mu_\alpha(\hhz(\eee))-\beta} <0 
$$
if $\mu_{\alpha}(\hhm(\eee))\neq\beta$ or, otherwise and similarly, $\infty\leq -(\mu_\alpha(\hhz(\eee))-\beta)^{-1} <0$. But this cannot happen unless one (and only one) among $\hhm(\eee)$ and $\hhz(\eee)$ vanishes. Say the latter does not, that is, 
$\mathcal{E}\cong \hh^{0}(\mathcal{E}) \in \mathcal{T}_\alpha^\beta$. 
If $\mathcal{E}$ is not $\mu_\alpha$-semistable, let $\mathcal{F}\in\csys(C)$ be a destabilizer, which we may further assume is $\mu_{\alpha}$-semistable. Then $\mu_\alpha(\mathcal{F})>\mu_\alpha(\eee)>\beta$. Thus Proposition \ref{prpelm}.(ii) yields $\mathcal{F}\in\mathcal{T}_\alpha^\beta\subset\csys_\alpha^\beta(C)$. Now
$$
\lim_{\gamma \to \infty} \frac{\mu^{\beta, \gamma}_\alpha(\mathcal{F})}{\gamma}=-\frac{1}{\mu_\alpha(\mathcal{F})-\beta}>-\frac{1}{\mu_\alpha(\mathcal{E})-\beta}=\lim_{\gamma \to \infty} \frac{\mu^{\beta, \gamma}_\alpha(\mathcal{E})}{\gamma}.
$$
In particular, there exists $\gamma_0\gg 0$ such that
$$
\frac{\mu^{\beta, \gamma_0}_\alpha(\mathcal{F})}{\gamma_0} > \frac{\mu^{\beta, \gamma_0}_\alpha(\mathcal{E})}{\gamma_0}
$$
which contradicts the $\mu^{\beta, \gamma_0}_\alpha$-semistability of $\mathcal{E}$. Therefore $\hh^{0}(\mathcal{E})$ is $\mu_\alpha$-semistable.

If $\eee=\hhm(\eee)[1]$ use a similar argument and Proposition \ref{prpelm}.(i) to conclude that $\hhm(\mathcal{E})$ is $\mu_\alpha$-semistable.
\end{proof}

Now for $\alpha, \beta \in \mathbb{Q}$, set 
\[
c=c_{\alpha}^{\beta,\gamma}:= {\rm min}\{\,\Im(Z^{\beta, \gamma}_\alpha(\mathcal{E}))\, |\, \mathcal{E} \in \csys_\alpha^\beta(C) \,\, \text{and} \,\, \Im(Z^{\beta, \gamma}_\alpha(\mathcal{E}))>0\}.\]

Note that this minimum exists since the image of $\Im(Z^{\beta, \gamma})$ is discrete if $\alpha,\beta\in\mathbb{Q}$ owing to Lemma \ref{discrete}.  In fact, its proof yields that if $\alpha=a_\alpha/b_\alpha$ and  $\beta=a_\beta/b_\beta$, then $c =N/(b_\alpha b_\beta)$ for some $N \in \mathbb{Z}_{>0}$. In particular, $c \in \mathbb{Q}_{>0}$.

\begin{Lemma}\label{minimalstable}
    Let $\alpha, \beta \in \mathbb{Q}$ and $\mathcal{E} \in \csys_\alpha^\beta(C)$ be such that  $\Im(Z^{\beta, \gamma}_\alpha(\mathcal{E}))=c$.   Then $\mathcal{E}$ is $\mu_\alpha^{\beta,\gamma}$-stable if and only if 
    $Hom (\mathcal{F},\mathcal{E})=0$ for all $\mathcal{F} \in \csys_\alpha^{\beta,\gamma}(C)$ with $\Im(Z^{\beta, \gamma}_\alpha(\mathcal{F}))=0$.
\end{Lemma}

\begin{proof}
First note necessity only requires $\mu_{\alpha}^{\beta,\gamma}(\eee)<\infty$ (which holds as $\Im(Z^{\beta, \gamma}_\alpha(\mathcal{E}))=c>0$) and semistability. Indeed, suppose so and consider $\mathcal{F}\in\csys_\alpha^\beta(C)$ with $\Im(Z^{\beta, \gamma}_\alpha(\mathcal{F}))=0$. Then $\mu_{\alpha}^{\beta,\gamma}(\mathcal{F})=\infty$ and, in particular, $\mathcal{F}$ is $\mu_\alpha^{\beta,\gamma}$-semistable. Then Proposition \ref{Homomorphism} yields $Hom (\mathcal{F},\mathcal{E})=0$.  To prove sufficiency, let  
\begin{equation}
\label{fit1}
0\to \mathcal{F} \to \mathcal{E} \to \mathcal{G} \to 0
\end{equation}    
be an exact triangle in $\csys_\alpha^\beta(C)$. Then $c=\Im(Z^{\beta, \gamma}_\alpha(\mathcal{E}))= \Im(Z^{\beta, \gamma}_\alpha(\mathcal{F})) + \Im(Z^{\beta, \gamma}_\alpha(\mathcal{G}))$. But, by Proposition \ref{WSC}, we have that $\Im(Z^{\beta, \gamma}_\alpha(\mathcal{F}))\geq 0$ and $\Im(Z^{\beta, \gamma}_\alpha(\mathcal{G}))\geq 0$. Hence either $\Im(Z^{\beta, \gamma}_\alpha(\mathcal{F}))=0$ or $\Im(Z^{\beta, \gamma}_\alpha(\mathcal{G}))=0$.
But $\Im(Z^{\beta, \gamma}_\alpha(\mathcal{F}))$ doesn't vanish by hypothesis. Thus $\Im(Z^{\beta, \gamma}_\alpha(\mathcal{G}))$ does. So $\mu_\alpha^{\beta,\gamma}(\mathcal{G})=\infty > \mu_\alpha^{\beta,\gamma}(\mathcal{E})$ and by the seesaw principle (\ref{fit1}) cannot destabilize $\mathcal{E}$. Therefore $\mathcal{E}$ is $\mu_\alpha^{\beta,\gamma}$-stable. 
\end{proof}

\begin{Lemma}
\label{semiminimal}
Let $\mathcal{E} \in \csys_\alpha^\beta(C)$. The following holds:
\begin{itemize}
    \item[(i)] if $\Im(Z^{\beta, \gamma}_\alpha(\mathcal{E})) = 0$, then $\mathcal{E}$ is $\mu_\alpha^{\beta,\gamma}$-semistable and $\mathcal{E} \in \textsf{B}$.
    \item[(ii)] if $\alpha,\beta\in\mathbb{Q}$ and $\Im(Z^{\beta, \gamma}_\alpha(\mathcal{E}))=c$, then $\mathcal{E}$ is $\mu_\alpha^{\beta, \gamma}$-semistable if and only if $\mathcal{E} \in \textsf{B}$.
\end{itemize} 
\end{Lemma}
\begin{proof}
Note that (i) easily follows from Corollary \ref{Im=0}, from the equality $\mu_{\alpha}^{\beta,\gamma}(\eee)=\infty$, and from the very definition of $\textsf{B}$.

To prove (ii), suppose $\mathcal{E}$ is $\mu_\alpha^{\beta,\gamma}$-semistable. Use the proof of necessity in Lemma \ref{minimalstable} to conclude that $Hom (\mathcal{F},\mathcal{E})=0$ for all $\mathcal{F} \in \csys_\alpha^{\beta,\gamma}(C)$ with $\Im(Z^{\beta, \gamma}_\alpha(\mathcal{F}))=0$. But this condition characterizes stability and does not depend on $\gamma$ as $\Im(Z^{\beta, \gamma}_\alpha)$ does not involve $\gamma$. Thus $\eee$ is $\mu_{\alpha}^{\beta,\gamma}$-stable for all $\gamma$. In particular, $\eee$ is  $\mu_\alpha^{\beta,\gamma}$-semistable for $\gamma \gg 0$ and hence $\mathcal{E} \in \textsf{B}$ by Lemma \ref{biggamma}. 

Conversely, assume $\mathcal{E} \in \textsf{B}$ 
and consider an exact triangle 
\begin{equation}
\label{fit1a}
0\to \mathcal{F} \to \mathcal{E} \to \mathcal{G} \to 0
\end{equation}    
in $\csys_\alpha^\beta(C)$. As $\Im(Z^{\beta, \gamma}_\alpha(\mathcal{E}))=c$, then either $\Im(Z^{\beta, \gamma}_\alpha(\mathcal{F}))=0$ or $\Im(Z^{\beta, \gamma}_\alpha(\mathcal{G}))=0$. If the latter holds then $\mu_{\alpha}^{\beta,\gamma}(\mathcal{F}) < \mu_{\alpha}^{\beta,\gamma}(\mathcal{E})$ by the seesaw principle as $\mu_{\alpha}^{\beta,\gamma}(\mathcal{G})=\infty$. So suppose $\Im(Z^{\beta, \gamma}_\alpha(\mathcal{F}))=0$. By Corollary \ref{Im=0}, $\mathcal{F} \cong \hh^{-1}(\mathcal{F})[1]$. Now $\eee\in\textsf{B}$ so either $\eee\cong\hhm(\eee)[1]$ or $\eee\cong\hhz(\eee)$. But the latter is precluded since, otherwise, \eqref{fit1a} yields $0\to\hhm(\mathcal{F})\to 0$ and hence $\mathcal{F}=0$. So $\eee\cong\hhm(\eee)[1]$ and \eqref{fit1a} yields an inclusion $\hhm(\mathcal{F})\subset\hhm(\eee)$. Now, by Corollary \ref{Im=0}, $\mu_\alpha(\hhm(\mathcal{F}))=\beta$. But $\hhm(\eee)\in\mathcal{F}_\alpha^\beta$ and is $\mu_\alpha$-semistable. Thus $\mu_\alpha(\hhz(\eee))=\beta$ too. This implies $\infty=\mu_{\alpha}^{\beta,\gamma}(\hhm(\eee))=\mu_{\alpha}^{\beta,\gamma}(\eee)$, which cannot happen as $\Im(Z^{\beta, \gamma}_\alpha(\mathcal{E}))=c>0$. So $\Im(Z^{\beta, \gamma}_\alpha(\mathcal{F}))$ doesn't vanish. Therefore, $\mathcal{E}$ is $\mu_\alpha^{\beta,\gamma}$-stable and, in particular, $\mu_\alpha^{\beta,\gamma}$-semistable. 
\end{proof}

\begin{Lemma}\label{subtor}
    Let $\eee \in \csys(C)$ be  $\mu_\alpha$-stable of rank $n>0$. Assume that $\mu_\alpha(\mathcal{E})\leq\beta$. Assume also that $\beta < \delta$ if $n\geq 2$, where $\delta=min\{\mu_\alpha(\mathcal{Q}) \, | \, \mathcal{E} \twoheadrightarrow \mathcal{Q}\neq \eee\}$. Then $\eee\in\fab$ and any proper subobject $\mathcal{F}\subset \mathcal{E}[1]$ in $\csys_\alpha^\beta(C)$ is in $\tab$ and fits into an exact sequence in $\csys(C)$ of the form 
    \begin{equation}\label{seq}
       0 \to \mathcal{E} \to \mathcal{F}' \to \mathcal{F} \to 0
    \end{equation}
    where $\mathcal{F}' \in \fab$.  
\end{Lemma}

\begin{proof}
As $\eee$ is $\mua$-stable and $\mua(\eee)\leq\beta$, it follows that $\eee\in\fab$. Thus $\eee[1]\in\csysab(C)$. Let $\fff$ be a proper subobject  of $\eee[1]$ in $\csysab(C)$, fitting into the exact sequence 
    \[
       0\to  \mathcal{F} \to \eee[1] \to \mathcal{G}\to 0
    \]
    in $\csys_\alpha^\beta(C)$. It yields the following exact sequence in $\csys(C)$ 
    \begin{equation}\label{fit}
    0 \to \hh^{-1}(\mathcal{F}) \stackrel{f}{\to}\mathcal{E} \to \hh^{-1}(\mathcal{G}) \to \hh^{0}(\mathcal{F})\to 0\to\hhz(\mathcal{G})\to 0.\end{equation}
      Thus $\hhz(\mathcal{G})=0$. Set $\mathcal{S}:={\rm coker}f$.  If $\mathcal{S}=0$, then $\hh^{-1}(\mathcal{G}) \cong \hh^{0}(\mathcal{F})$ which implies $\hhm(\mathcal{G})=0$, and hence $\mathcal{G}=0$ since its cohomologies vanish. But this is impossible as $\fff$ is a proper subobject of $\mathcal{E}[1]$.  Therefore $\mathcal{S}\neq 0$. 

      If $n=1$, we have $\mathcal{S}\subset \hh^{-1}(\mathcal{G})\in \mathcal{F}_\alpha^\beta$. Thus $\rank(\mathcal{S})\neq 0$. So $\rank(\hh^{-1}(\mathcal{F})) =0$ since $\rank(\mathcal{E})=1$. But $\hh^{-1}(\mathcal{F})\in \mathcal{F}_\alpha^\beta$, and so $\hh^{-1}(\mathcal{F})=0$. Thus $\fff=\hh^{0}(\mathcal{F})\in \mathcal{T}_\alpha^\beta$ as desired.
      
    In $n\geq 2$, assume $\mathcal{S}$ is a proper quotient of $\eee$. Then $\mu_\alpha(\eee) \leq \beta < 
    \delta \leq \mu_\alpha(\mathcal{S})$. On the other hand, $\mu_\alpha(\mathcal{S}) \leq \beta$ because $\mathcal{S} \hookrightarrow \hh^{-1}(\mathcal{G}) \in \mathcal{F}_\alpha^\beta$ which yields a contradiction.  Therefore $\mathcal{S}=\eee$ and hence $\hh^{-1}(\mathcal{F})= 0$. Thus $\mathcal{F}=\hh^0(\mathcal{F})\in\tab$ as desired. Define $\fff':=\hhm(\mathcal{G})$ and \eqref{fit} yields  \eqref{seq}.
\end{proof}

\begin{Proposition}
\label{semisimple}
Let $\mathcal{E}\in\csys(C)$ be $\mu_\alpha$-stable with $\mu_\alpha(\mathcal{E})=\beta$. Then $\mathcal{E}[1]$ is minimal in $\csysab(C)$. 
 \end{Proposition}

 \begin{proof}
By Proposition \ref{prpelm} (or even Lemma \ref{subtor}), we have that $\eee\in\fab$. So let  $\mathcal{F}$ be a proper subobject of $\eee[1]$. By Lemma \ref{subtor}, $\fff\in\tab$ and fits into the exact sequence 
      \[0 \to \mathcal{E} \to \mathcal{F}' \to \mathcal{F} \to 0\]
      in $\csys(C)$ with $\fff'\in\fab$. But the seesaw principle, along with the definitions of $\fab$ and $\tab$, yields:
$$ \mu_\alpha(\mathcal{E})=\beta\geq  \mu_\alpha(\mathcal{F}') \geq \mu_\alpha(\mathcal{F})>\beta $$
which is a contradiction. Thus $\mathcal{F}$ is not proper, that is, $\mathcal{E}[1]$ is minimal.
\end{proof}

 We can now formulate our main result of this section.

  \begin{Theorem}\label{stabpos}
     Let $\alpha,  \gamma \in \mathbb{Q}$ . Let $\eee \in \csys(C)$ be  $\mu_\alpha$-stable of type $(n,d,k)$, with $n>0$ and $-\gamma n<d$ if $k=0$. Assume  $\mu_\alpha(\mathcal{E}) <\beta$, $\beta \in \mathbb{Q}$. Assume also that $\beta < \delta_\eee$ if $n\geq 2$, where $\delta_\eee=min\{\mu_\alpha(\mathcal{Q}) \, | \, \mathcal{E} \twoheadrightarrow \mathcal{Q}\neq \eee\}$. Then; 
     \begin{itemize}
         \item [(i)] If $\eee=(E,V)$ is not complete, then $\mathcal{E}[1]$ is $\mu_\alpha^{\beta,\gamma}$-stable  in $\csys_\alpha^\beta(C)$ for every $ \gamma \gg 0$  and $\mu_\alpha^{\beta,\gamma}$-unstable for $ \gamma < \gamma_0:=(\beta n-d(\alpha+1))/\alpha n$. Also, if $|\eee| \in \fab$, then $\gamma_0$ is a pseudo-wall for $(n,d,k)$. Moreover, if $|\eee|$ is $\mu_\alpha$-stable and $\beta <\delta_{|\eee|}$, then  $\gamma_0$ is an actual-wall for $(n,d,k)$
         \item [(ii)] If $\eee$ is  complete, then $\mathcal{E}[1]$ is $\mu_\alpha^{\beta,\gamma}$-stable for every $\gamma$ .
     \end{itemize}
\end{Theorem}

\begin{proof}
 First note that $\eee \in \mathscr{F}_\alpha^\beta$. Let $\fff$ be a  subobjects of $\eee[1]$ in $\csysab(C)$.    From Lemma \ref{subtor}, we have $\mathcal{F} \in \tab$ and fits into an  exact sequence 
      \[0 \to \mathcal{E} \to \mathcal{F}' \to \mathcal{F} \to 0\]
      in $\csys(C)$ with  $\mathcal{F}' \in \fab$
     which induces the exact sequence 
     \[0\to \mathcal{F} \to \mathcal{E}[1] \to \mathcal{F}'[1] \to 0\]
     in $\csys_\alpha^\beta(C)$.  Now say $\fff$ is of type $(n',d',k')$. As $\mua(\fff)>\beta\geq 0$, it follows that $d' \geq 0$ if $k'=0$ and also $\mu_\alpha^{\beta, \gamma}(\fff)\neq \infty$. If $\fff$ is either pure $n'>0$, or not pure and injective, or $n'=0$ and $d'\geq k'$, then $\mu_\alpha^{\beta,\gamma}(\fff)\leq 0$ by Proposition \ref{TorHNF}(ii). If $\fff$ is not pure and not injective with $n>0$, then $k'\leq d'+n'+u-t$ for some $t, u \in \mathbb{Z}$. Hence, there exists $\gamma \gg 0$ such that $k'\leq d'+n'+u-t \leq d'+\gamma n'$. Thus, $\mu_\alpha^{\beta,\gamma}(\fff)\leq 0$. On the other hand, $\eee$ is $\mua$-stable with $\mu_a(\eee) \leq \beta$. So Proposition \ref{slopes}(iii) and our assumption on $d$ and $k$, yield that $\mu_\alpha^{\beta,\gamma}(\eee)>0$ and we are done.
     
     So assume $n'=0$ and $0\leq d'<k'$.  For the $\mu_\alpha^{\beta,\gamma}$-stability of $\eee[1]$, we need
\begin{equation}
\label{equin111}
\mu_\alpha^{\beta,\gamma}(\fff)=\frac{k'-d'}{d'+\alpha k'} <  \frac{d+\gamma n-k}{\beta n-(d+\alpha k)}=\mu_\alpha^{\beta,\gamma}(\eee)=\mu_\alpha^{\beta,\gamma}(\eee[1])
\end{equation}
But we have that
\begin{equation}
\label{equin211}
\frac{k'-d'}{d'+\alpha k'}=\frac{1-d'/k'}{\alpha+d'/k'} \leq \frac{1}{\alpha}
\end{equation}
Combining \eqref{equin111} and \eqref{equin211}, it suffices for us to have
\begin{equation}
\label{equin31}
\frac{\beta n-d(\alpha+1)}{\alpha n} < \gamma 
\end{equation}
Combining the two cases, we have that $\eee[1]$ is $\mu_\alpha^{\beta,\gamma}$-stable for $\gamma \gg 0$ as we desired.

To prove (i), write $\eee=(E,V)$ and assume $\eee$ is not complete. Then we can consider the following exact sequence in $\csys(C)$,
$$ 0 \to \mathcal{E} \to |\eee| \to (0,H^0(E)/V) \to 0. $$
whose last term does not vanish. As $|\eee| \in \fab$, the above sequence yields the following  exact sequence in $\csys_{\alpha}^\beta(C)$,
\begin{equation} \label{wallnoncomplete}
0\to (0,H^0(E)/V) \to \mathcal{E}[1] \to |\eee|[1]  \to 0 . 
\end{equation}
Note that 
$$ \mu_\alpha^{\beta,\gamma}(0,H^0(E)/V) = \dfrac{1}{\alpha} ~~,~~ \mu_\alpha^{\beta,\gamma}(\mathcal{E}[1]) = \dfrac{-(d+\gamma n -k)}{d+\alpha k-\beta n}. $$
Therefore, $1/\alpha > -(d+\gamma n-k)/(d+\alpha k-\beta n)$ since $\gamma <(\beta n-d(\alpha+1)/\alpha n$. 
Thus $\eee[1]$ is $\mu_\alpha^{\beta,\gamma}$-unstable and the surface $\gamma_0 = (\beta n-d(\alpha+1))/\alpha n,$ is a pseudo wall for $(n,d,k)$. To prove that $\gamma_0$ is an actual wall, we first prove item (ii).

To prove (ii),  write $\mathcal{F}=(F,U,\phi)$ and $\mathcal{F}'=(F',U',\phi')$, this yields the following commutative diagram
$$ \xymatrix{
0 \ar[r] & V \ar[r] \ar[d]^{\varphi} & U' \ar[r] \ar[d]^{\phi'} & U \ar[r] \ar[d]^{\phi} & 0 \ar[d] \\ 
0 \ar[r] & H^0(E) \ar[r] & H^0(F') \ar[r] & H^0(F) \ar[r] & H^1(E)
}$$
By the Four-lemma, the morphism $\phi$ is injective since $\varphi$ is an isomorphism. Therefore, $\mathcal{F}$ is injective. Say $\mathcal{F}$ is of type $(n',d',k')$. Hence  $k'\le d'+n'$ by Proposition \ref{TorHNF}.(ii)(b). From the above it follows that $\mu_\alpha^{\beta,\gamma}(\mathcal{F}) <0 <\mu_\alpha^{\beta,\gamma}(\eee[1])$ as desired.

Now we proceed to show that $\gamma_0$ is an actual wall in item (i). From item (ii), we have that in the exact sequence (\ref{wallnoncomplete}), the object $|\eee|[1]$ is $\mu_\alpha^{\beta,\gamma}$-stable for any $\gamma$, and also, from Theorem \ref{geometric} we have that $(0,H^0(E)/V)$ is $\mu_\alpha^{\beta,\gamma}$-semistable.  Moreover for $\gamma=\gamma_0$, we have \[\mu_\alpha^{\beta,\gamma}(0,H^0(E)/V) = \mu_\alpha^{\beta,\gamma}(\eee[1]).\] Hence, it follows easily that $\eee[1]$ is $\mu_\alpha^{\beta,\gamma}$-semistable. Therefore,  for $\gamma > \gamma_0$ we have that $\eee[1]$ is $\mu_\alpha^{\beta,\gamma}$-stable and we can conclude that $\gamma_0$ is an actual wall and also for $\gamma < \gamma_0$  which is the desired conclusion.  
\end{proof}




\begin{center} References \end{center}
\begin{biblist}




\bib{ACGH}{book} {
    AUTHOR = {Arbarello, E. and Cornalba, M. and Griffiths, P. A. and Harris, J.},
     TITLE = {Geometry of algebraic curves. {V}ol. {I}},
    SERIES = {Grundlehren der mathematischen Wissenschaften [Fundamental
              Principles of Mathematical Sciences]},
    VOLUME = {267},
 PUBLISHER = {Springer-Verlag, New York},
      YEAR = {1985},
     PAGES = {xvi+386},
      ISBN = {0-387-90997-4},
       DOI = {10.1007/978-1-4757-5323-3},
       URL = {https://doi.org/10.1007/978-1-4757-5323-3},
}

\bib{ABCH}{article}{
   author={Arcara, Daniele},
   author={Bertram, Aaron},
   author={Coskun, Izzet},
   author={Huizenga, Jack},
   title={The minimal model program for the Hilbert scheme of points on
   $\mathbb{P}^2$ and Bridgeland stability},
   journal={Adv. Math.},
   volume={235},
   date={2013},
   pages={580--626},
   issn={0001-8708},
   review={\MR{3010070}},
   doi={10.1016/j.aim.2012.11.018},
}

\bib{Bayer}{article}{
AUTHOR = {A. Bayer.},
TITLE = {A tour to stability conditions on derived categories. http://www.maths.ed.ac.uk/~abayer/dc-lecture-notes.pdf, 2011.},
}

\bib{BM11}{article}{
author = {Arend Bayer and Emanuele Macr{\`i}},
title = {{The space of stability conditions on the local projective plane}},
volume = {160},
journal = {Duke Mathematical Journal},
number = {2},
publisher = {Duke University Press},
pages = {263 -- 322},
year = {2011},
doi = {10.1215/00127094-1444249},
URL = {https://doi.org/10.1215/00127094-1444249}
}

\bib{BM14a}{article}{
	Author = {Bayer, Arend and Macr{\`\i}, Emanuele},
	Journal = {Inventiones mathematicae},
	Number = {3},
	Pages = {505--590},
	Title = {MMP for moduli of sheaves on K3s via wall-crossing: nef and movable cones, Lagrangian fibrations},
	Volume = {198},
	Year = {2014},
	Issn = {1432-1297},
	DOI = {10.1007/s00222-014-0501-8},
	URL = {https://doi.org/10.1007/s00222-014-0501-8}}

\bib{BMT}{article}{
   author={Bayer, Arend},
   author={Macr\`i, Emanuele},
   author={Toda, Yukinobu},
   title={Bridgeland stability conditions on threefolds I:
   Bogomolov-Gieseker type inequalities},
   journal={J. Algebraic Geom.},
   volume={23},
   date={2014},
   number={1},
   pages={117--163},
   issn={1056-3911},
   review={\MR{3121850}},
   doi={10.1090/S1056-3911-2013-00617-7},
}

\bib{BMS}{article}{
  author  = {Bayer, Arend and Macr\`i, Emanuele and Stellari, Paolo},
  title   = {The space of stability conditions on abelian threefolds, and on some {Calabi-Yau} threefolds},
  journal = {Inventiones Mathematicae},
  volume  = {214},
  number  = {1},
  pages   = {303--348},
  year    = {2018},
  doi     = {10.1007/s00222-018-0797-0},
}

\bib{BGPMN}{article}{
   author={Bradlow, S. B.},
   author={Garc\'ia-Prada, O.},
   author={Mu\~noz, V.},
   author={Newstead, P. E.},
   title={Coherent systems and Brill-Noether theory},
   journal={Internat. J. Math.},
   volume={14},
   date={2003},
   number={7},
   pages={683--733},
   issn={0129-167X},
   review={\MR{2001263}},
   doi={10.1142/S0129167X03002009},
}

\bib{BDG}{article}{
   author={Burban, I.},
   author={Drozd, Yu.},
   author={Greuel, G.-M.},
   title={Vector bundles on singular projective curves},
   conference={
      title={Applications of algebraic geometry to coding theory, physics
      and computation},
      address={Eilat},
      date={2001},
   },
   book={
      series={NATO Sci. Ser. II Math. Phys. Chem.},
      volume={36},
      publisher={Kluwer Acad. Publ., Dordrecht},
   },
   isbn={1-4020-0004-9},
   date={2001},
   pages={1--15},
   review={\MR{1866891}},
}

\bib{BGMN2}{article}{
   author={Bradlow, Steven},
   author={Daskalopoulos, Georgios D.},
   author={Garc\'ia-Prada, Oscar},
   author={Wentworth, Richard},
   title={Stable augmented bundles over Riemann surfaces},
   conference={
      title={Vector bundles in algebraic geometry},
      address={Durham},
      date={1993},
   },
   book={
      series={London Math. Soc. Lecture Note Ser.},
      volume={208},
      publisher={Cambridge Univ. Press, Cambridge},
   },
   isbn={0-521-49878-3},
   date={1995},
   pages={15--67},
   review={\MR{1338412}},
   doi={10.1017/CBO9780511569319.003},
}


\bib{BGMN1}{article}{
   author={Bradlow, Steven B.},
   author={Garc\'ia-Prada, Oscar},
   title={An application of coherent systems to a Brill-Noether problem},
   journal={J. Reine Angew. Math.},
   volume={551},
   date={2002},
   pages={123--143},
   issn={0075-4102},
   review={\MR{1932176}},
   doi={10.1515/crll.2002.079},
}

\bib{BGN}{article} {
    AUTHOR = {Brambila-Paz L.},
    AUTHOR = {Grzegorczyk, I.},
    AUTHOR = {Newstead, P. E.},
     TITLE = {Geography of {B}rill-{N}oether loci for small slopes},
   JOURNAL = {J. Algebraic Geom.},
  JOURNAL = {Journal of Algebraic Geometry},
    VOLUME = {6},
      YEAR = {1997},
    NUMBER = {4},
     PAGES = {645--669},
      ISSN = {1056-3911,1534-7486},
}

\bib{BramLange}{article}{
   author={Brambila-Paz, L.},
   author={Lange, Herbert},
   title={A stratification of the moduli space of vector bundles on curves},
   note={Dedicated to Martin Kneser on the occasion of his 70th birthday},
   journal={J. Reine Angew. Math.},
   volume={494},
   date={1998},
   pages={173--187},
   issn={0075-4102},
   review={\MR{1604403}},
   doi={10.1515/crll.1998.005},
}

\bib{Bridgeland}{article}{
   author={Bridgeland, Tom},
   title={Stability conditions on triangulated categories},
   journal={Ann. of Math. (2)},
   volume={166},
   date={2007},
   number={2},
   pages={317--345},
   issn={0003-486X},
   review={\MR{2373143}},
   doi={10.4007/annals.2007.166.317},
}

\bib{CR}{article}{
   author={Cotterill, Ethan},
   author={Martins, Renato Vidal},
   title={Towards Brill-Noether theory for cuspidal curves},
   journal={Mat. Contemp.},
   volume={60},
   date={2024},
   pages={31--48},
   issn={0103-9059},
   review={\MR{4729077}},
   doi={10.21711/231766362024/rmc603},
}

\bib{Bridgeland1}{article}{
    AUTHOR = {Bridgeland, Tom},
     TITLE = {Stability conditions on {$K3$} surfaces},
   JOURNAL = {Duke Math. J.},
  JOURNAL = {Duke Mathematical Journal},
    VOLUME = {141},
      YEAR = {2008},
    NUMBER = {2},
     PAGES = {241--291},
      ISSN = {0012-7094,1547-7398},
       DOI = {10.1215/S0012-7094-08-14122-5},
       URL = {https://doi.org/10.1215/S0012-7094-08-14122-5},
}



\bib{EKS}{article}{
   author={Eisenbud, David},
   author={Koh, Jee},
   author={Stillman, Michael},
   title={Determinantal equations for curves of high degree},
   journal={Amer. J. Math.},
   volume={110},
   date={1988},
   number={3},
   pages={513--539},
   issn={0002-9327},
   doi={10.2307/2374621},
}

\bib{LGMS}{article}{
author={Feital, L.},
author={Galdino, N.},
    author = {Martins, R. V.},
    author = {Souza, A.},
     title = {On Clifford dimension for singular curves, \emph{}},
      }
  \bib{FL}{article}{
   author={Fulton, W.},
   author={Lazarsfeld, R.},
   title={On the connectedness of degeneracy loci and special divisors},
   journal={Acta Math.},
   volume={146},
   date={1981},
   number={3-4},
   pages={271--283},
   issn={0001-5962},
   review={\MR{611386}},
   doi={10.1007/BF02392466},
}
\bib{FN}{article}{
      title={Derived category of coherent systems on curves and stability conditions}, 
      author={Feyzbakhsh, Soheyla},
      author={Novik, Aliaksandra},
      year={2025},
      eprint={2511.01601},
      url={https://arxiv.org/abs/2511.01601}, 
}

\bib{GH}{article}{
   author={Griffiths, Phillip},
   author={Harris, Joseph},
   title={On the variety of special linear systems on a general algebraic
   curve},
   journal={Duke Math. J.},
   volume={47},
   date={1980},
   number={1},
   pages={233--272},
   issn={0012-7094},
   review={\MR{0563378}},
}





\bib{H}{book}{
   author={Hartshorne, Robin},
   title={Algebraic geometry},
   note={Graduate Texts in Mathematics, No. 52},
   publisher={Springer-Verlag, New York-Heidelberg},
   date={1977},
   pages={xvi+496},
   isbn={0-387-90244-9},
   review={\MR{0463157}},
}






\bib{He}{article}{
   author={He, Min},
   title={Espaces de modules de syst\`emes coh\'erents},
   language={French},
   journal={Internat. J. Math.},
   volume={9},
   date={1998},
   number={5},
   pages={545--598},
   issn={0129-167X},
   review={\MR{1644040}},
   doi={10.1142/S0129167X98000257},
}

\bib{HL}{book}{
   author={Huybrechts, Daniel},
   author={Lehn, Manfred},
   title={The geometry of moduli spaces of sheaves},
   series={Cambridge Mathematical Library},
   edition={2},
   publisher={Cambridge University Press, Cambridge},
   date={2010},
   pages={xviii+325},
   isbn={978-0-521-13420-0},
   review={\MR{2665168}},
   doi={10.1017/CBO9780511711985},
}

\bib{J}{book}{
   author={Jacobson, Nathan},
   title={Basic algebra. I},
   edition={2},
   publisher={W. H. Freeman and Company, New York},
   date={1985},
   pages={xviii+499},
   isbn={0-7167-1480-9},
   review={\MR{0780184}},
}

\bib{JLMM}{article}{
      title={Higher rank DT/PT wall-crossing in Bridgeland stability}, 
      author={Jardim, Marcos},
      author={Lo, Jason},
      author={Maciocia, Antony},
      author={Martinez, Cristian},
      year={2025},
      eprint={2503.20008},
      url={https://arxiv.org/abs/2503.20008}, 
}

\bib{JM}{article}{
   author={Jardim, Marcos},
   author={Maciocia, Antony},
   title={Walls and asymptotics for Bridgeland stability conditions on
   threefolds},
   journal={\'Epijournal G\'eom. Alg\'ebrique},
   volume={6},
   date={2022},
   pages={Art. 22, 61},
   review={\MR{4526270}},
}

\bib{Ke}{article}{
   author={Kempf, G.},
   title={schubert methods with an application to algebraic curves},
   journal={Public. Math. Centrum.},
   volume={},
   date={1971},
   number={},
   pages={},
   issn={},
   review={Amsterdam},
   doi={},
}


\bib{KN}{article}{
   author={King, A. D.},
   author={Newstead, P. E.},
   title={Moduli of Brill--Noether pairs on algebraic curves},
   journal={Internat. J. Math.},
   volume={6},
   date={1995},
   number={5},
   pages={733--748},
   issn={0129-167X},
   doi={10.1142/S0129167X95000316},
}



\bib{KL72}{article}{
   author={Kleiman, Steven L.},
   author={Laksov, Dan},
   title={On the existence of special divisors},
   journal={Amer. J. Math.},
   volume={94},
   date={1972},
   pages={431--436},
   issn={0002-9327},
   review={\MR{0323792}},
   doi={10.2307/2374630},
}

\bib{KL74}{article}{
   author={Kleiman, Steven L.},
   author={Laksov, Dan},
   title={Another proof of the existence of special divisors},
   journal={Acta Math.},
   volume={132},
   date={1974},
   pages={163--176},
   issn={0001-5962},
   review={\MR{357398}},
   doi={10.1007/BF02392112},
}


\bib{KM}{article}{
   author={Kleiman, Steven Lawrence},
   author={Martins, Renato Vidal},
   title={The canonical model of a singular curve},
   journal={Geom. Dedicata},
   volume={139},
   date={2009},
   pages={139--166},
   issn={0046-5755},
   review={\MR{2481842}},
   doi={10.1007/s10711-008-9331-4},
}

\bibitem{KM2} S. L. Kleiman and R. V. Martins, {\it Gonality and Scrolls of an Integral Curve}, work in progress.



\bib{LN1}{article}{
    AUTHOR = {Lange, Herbert},
    author= { Newstead, Peter E.},
     TITLE = {Clifford indices for vector bundles on curves},
 BOOKTITLE = {Affine flag manifolds and principal bundles},
    SERIES = {Trends Math.},
     PAGES = {165--202},
 PUBLISHER = {Birkh\"auser/Springer Basel AG, Basel},
      YEAR = {2010},
      ISBN = {978-3-0346-0287-7},
       DOI = {10.1007/978-3-0346-0288-4\_6},
       URL = {https://doi.org/10.1007/978-3-0346-0288-4_6},
}

\bib{LNg0}{article}{
    AUTHOR = {Lange, H.},
    AUTHOR = {Newstead, P. E.},
     TITLE = {Coherent systems of genus 0. {II}. {E}xistence results for
              {$k\geq 3$}},
   JOURNAL = {Internat. J. Math.},
  JOURNAL = {International Journal of Mathematics},
    VOLUME = {18},
      YEAR = {2007},
    NUMBER = {4},
     PAGES = {363--393},
      ISSN = {0129-167X,1793-6519},
       DOI = {10.1142/S0129167X07004072},
       URL = {https://doi.org/10.1142/S0129167X07004072},
}

\bib{Li}{article}{
  author  = {Li, Chunyi},
  title   = {On stability conditions for the quintic threefold},
  journal = {Mathematische Annalen},
  volume  = {376},
  number  = {3--4},
  pages   = {1615--1659},
  year    = {2020},
  doi     = {10.1007/s00222-019-00888-z},
}


\bib{Lu}{book}{
   author={L\"ubke, Martin},
   author={Teleman, Andrei},
   title={The Kobayashi-Hitchin correspondence},
   publisher={World Scientific Publishing Co., Inc., River Edge, NJ},
   date={1995},
   pages={x+254},
   isbn={981-02-2168-1},
   review={\MR{1370660}},
   doi={10.1142/2660},
}

\bib{JNO}{article}{
author={Neulaender, G.},
    author = {Oliveira, E.},
     title = {Stability conditions on abelian comma categories, \emph{Preprint arXiv:2510.25450}},
       YEAR = {2025},
      }

\bib{LN}{article}{
    author = {Lange, H.},
    author = {Newstead, P. E.},
     TITLE = {Clifford's theorem for coherent systems},
   JOURNAL = {Arch. Math. (Basel)},
  JOURNAL = {Archiv der Mathematik},
    VOLUME = {90},
      YEAR = {2008},
    NUMBER = {3},
     PAGES = {209--216},
      ISSN = {0003-889X,1420-8938},
       DOI = {10.1007/s00013-007-2534-3},
       URL = {https://doi.org/10.1007/s00013-007-2534-3},
}


\bib{LP}{article}{
   author={Le Potier, Joseph},
   title={Syst\`emes coh\'erents et structures de niveau},
   language={French, with English and French summaries},
   journal={Ast\'erisque},
   number={214},
   date={1993},
   pages={143},
   issn={0303-1179},
}



\bib{Macri-Schmidt}{article}{
   author={Macr\`i, Emanuele},
   author={Schmidt, Benjamin},
   title={Lectures on Bridgeland stability},
   conference={
      title={Moduli of curves},
   },
   book={
      series={Lect. Notes Unione Mat. Ital.},
      volume={21},
      publisher={Springer, Cham},
   },
   isbn={978-3-319-59485-9},
   isbn={978-3-319-59486-6},
   date={2017},
   pages={139--211},
}


\bib{Mer}{article}{
   author={Mercat, Vincent},
   title={Le probl\`eme de Brill-Noether pour des fibr\'es stables de petite
   pente},
   language={French},
   journal={J. Reine Angew. Math.},
   volume={506},
   date={1999},
   pages={1--41},
   issn={0075-4102},
   review={\MR{1665673}},
   doi={10.1515/crll.1999.005},
}

\bib{Mer2}{article}{
   author={Mercat, Vincent},
   title={Clifford's theorem and higher rank vector bundles},
   journal={Internat. J. Math.},
   volume={13},
   date={2002},
   number={7},
   pages={785--796},
   issn={0129-167X},
   review={\MR{1921510}},
   doi={10.1142/S0129167X02001484},
}


\bib{NarRam}{article}{
   author={Narasimhan, M. S.},
   author={Ramanan, S.},
   title={Geometry of Hecke cycles. I},
   conference={
      title={C. P. Ramanujam---a tribute},
   },
   book={
      series={Tata Inst. Fundam. Res. Stud. Math.},
      volume={8},
      publisher={Springer, Berlin-New York},
   },
   isbn={3-540-08770-2},
   date={1978},
   pages={291--345},
   review={\MR{0541029}},
}



\bib{Ne}{article}{
    AUTHOR = {Newstead, P. E.},
     TITLE = {Existence of {$\alpha$}-stable coherent systems on algebraic
              curves},
 BOOKTITLE = {Grassmannians, moduli spaces and vector bundles},
    SERIES = {Clay Math. Proc.},
    VOLUME = {14},
     PAGES = {121--139},
 PUBLISHER = {Amer. Math. Soc., Providence, RI},
      YEAR = {2011},
      ISBN = {978-0-8218-5205-7}
}

\bib{NMo}{article}{
    AUTHOR = {King, A. D. and Newstead, P. E.},
     TITLE = {Moduli of {B}rill-{N}oether pairs on algebraic curves},
   JOURNAL = {Internat. J. Math.},
  JOURNAL = {International Journal of Mathematics},
    VOLUME = {6},
      YEAR = {1995},
    NUMBER = {5},
     PAGES = {733--748},
      ISSN = {0129-167X,1793-6519},
       DOI = {10.1142/S0129167X95000316},
       URL = {https://doi.org/10.1142/S0129167X95000316},
}



\bib{New}{article}{
    AUTHOR = {Newstead, P. E.},
     TITLE = {Higher rank {B}rill-{N}oether theory and coherent systems open questions},
   JOURNAL = {Proyecciones},
  JOURNAL = {Proyecciones. Journal of Mathematics},
    VOLUME = {41},
      YEAR = {2022},
    NUMBER = {2},
     PAGES = {449--480},
      ISSN = {0716-0917,0717-6279},
}

\bib{PT}{article}{
    AUTHOR = {Piyaratne, Dulip},
     AUTHOR = {Toda, Yukinobu},
     TITLE = {Moduli of {B}ridgeland semistable objects on 3-folds and
              {D}onaldson-{T}homas invariants},
   JOURNAL = {J. Reine Angew. Math.},
  JOURNAL = {Journal f\"ur die Reine und Angewandte Mathematik. [Crelle's
              Journal]},
    VOLUME = {747},
      YEAR = {2019},
     PAGES = {175--219},
      ISSN = {0075-4102,1435-5345},
       DOI = {10.1515/crelle-2016-0006},
       URL = {https://doi.org/10.1515/crelle-2016-0006},
}



\bib{RHR}{article}{
AUTHOR = {Martínez-Romero, Eva},
AUTHOR = {Rincón-Hidalgo, Aejandra},
AUTHOR = {Rüffer, Arne},
      title={Bridgeland stability conditions on the category of holomorphic triples over curves}, 
      year={2020},
      eprint={1905.04240},
      url={https://arxiv.org/abs/1905.04240}, 
}

\bib{RaVi}{article}{
   author={Raghavendra, Nyshadham},
   author={Vishwanath, Periyapatna A.},
   title={Moduli of pairs and generalized theta divisors},
   journal={Tohoku Math. J. (2)},
   volume={46},
   date={1994},
   number={3},
   pages={321--340},
   issn={0040-8735},
   review={\MR{1289182}},
   doi={10.2748/tmj/1178225715},
}

\bib{RSt}{article}{
   author={Rosa, Renata},
   author={St\"{o}hr, Karl-Otto},
   title={Trigonal Gorenstein curves},
   journal={J. Pure Appl. Algebra},
   volume={174},
   date={2002},
   number={2},
   pages={187--205},
   issn={0022-4049},
   review={\MR{1921820}},
   doi={10.1016/S0022-4049(02)00122-6},
}

\bib{Th}{article}{
   author={Thaddeus, Michael},
   title={Stable pairs, linear systems and the Verlinde formula},
   journal={Invent. Math.},
   volume={117},
   date={1994},
   number={2},
   pages={317--353},
   issn={0020-9910},
   review={\MR{1273268}},
   doi={10.1007/BF01232244},
}

\bib{Ty}{article}{
   author={Tyurin, A. N.},
   title={Symplectic structures on the moduli spaces of vector bundles on
   algebraic surfaces with $p_g>0$},
   language={Russian},
   journal={Izv. Akad. Nauk SSSR Ser. Mat.},
   volume={52},
   date={1988},
   number={4},
   pages={813--852, 896},
   issn={0373-2436},
   translation={
      journal={Math. USSR-Izv.},
      volume={33},
      date={1989},
      number={1},
      pages={139--177},
      issn={0025-5726},
   },
   review={\MR{0966986}},
   doi={10.1070/IM1989v033n01ABEH000818},
}

\end{biblist}
\end{document}